\documentclass[12pt]{amsart}

\usepackage{amssymb}
\batchmode

\begin{document}

\def\dbl{[\hskip -1pt[}
\def\dbr{]\hskip -1pt]}
\title{Algebraic approximation in CR geometry}
\author{Nordine Mir}
\address{Universit\'e de Rouen, CNRS, Laboratoire de Math\'ematiques Rapha\"el Salem, UMR 6085 CNRS, Avenue de
l'Universit\'e, B.P. 12, 76801 Saint Etienne du Rouvray, France}
\email{Nordine.Mir@univ-rouen.fr}
\thanks{The author was partially supported by the French National Agency for Research (ANR), projects ANR-10-BLAN-0102 and  ANR-09-BLAN-0422.}



\def\Label#1{\label{#1}}

\def\1#1{\ov{#1}}
\def\2#1{\widetilde{#1}}
\def\6#1{\mathcal{#1}}
\def\4#1{\mathbb{#1}}
\def\3#1{\widehat{#1}}
\def\K{{\4K}}
\def\LL{{\4L}}

\def \MM{{\4M}}
\def \S{{\4S}^{2N'-1}}

\def \B{{\4B}^{2N'-1}}

\def \H{{\4H}^{2l-1}}

\def \F{{\4H}^{2N'-1}}

\def \LL{{\4L}}

\def\Re{{\sf Re}\,}
\def\Im{{\sf Im}\,}
\def\id{{\sf id}\,}

\def\s{s}
\def\k{\kappa}
\def\ov{\overline}
\def\span{\text{\rm span}}
\def\ad{\text{\rm ad }}
\def\tr{\text{\rm tr}}
\def\xo {{x_0}}
\def\Rk{\text{\rm Rk\,}}
\def\sg{\sigma}
\def \emxy{E_{(M,M')}(X,Y)}
\def \semxy{\scrE_{(M,M')}(X,Y)}
\def \jkxy {J^k(X,Y)}
\def \gkxy {G^k(X,Y)}
\def \exy {E(X,Y)}
\def \sexy{\scrE(X,Y)}
\def \hn {holomorphically nondegenerate}
\def\hyp{hypersurface}
\def\prt#1{{\partial \over\partial #1}}
\def\det{{\text{\rm det}}}
\def\wob{{w\over B(z)}}
\def\co{\chi_1}
\def\po{p_0}
\def\fb {\bar f}
\def\gb {\bar g}
\def\Fb {\ov F}
\def\Gb {\ov G}
\def\Hb {\ov H}
\def\zb {\bar z}
\def\wb {\bar w}
\def \qb {\bar Q}
\def \t {\tau}
\def\z{\chi}
\def\w{\tau}
\def\Z{\zeta}
\def\phi{\varphi}
\def\eps{\epsilon}

\def \T {\theta}
\def \Th {\Theta}
\def \L {\Lambda}
\def\b {\beta}
\def\a {\alpha}
\def\o {\omega}
\def\l {\lambda}

\def \im{\text{\rm Im }}
\def \re{\text{\rm Re }}
\def \Char{\text{\rm Char }}
\def \supp{\text{\rm supp }}
\def \codim{\text{\rm codim }}
\def \Ht{\text{\rm ht }}
\def \Dt{\text{\rm dt }}
\def \hO{\widehat{\mathcal O}}
\def \cl{\text{\rm cl }}
\def \bS{\mathbb S}
\def \bK{\mathbb K}
\def \bD{\mathbb D}
\def \bC{\mathbb C}
\def \bL{\mathbb L}
\def \bZ{\mathbb Z}
\def \bN{\mathbb N}
\def \scrF{\mathcal F}
\def \scrK{\mathcal K}
\def \mc #1 {\mathcal {#1}}
\def \scrM{\mathcal M}
\def \cR{\mathcal R}
\def \scrJ{\mathcal J}
\def \scrA{\mathcal A}
\def \scrO{\mathcal O}
\def \scrV{\mathcal V}
\def \scrL{\mathcal L}
\def \scrE{\mathcal E}
\def \hol{\text{\rm hol}}
\def \aut{\text{\rm aut}}
\def \Aut{\text{\rm Aut}}
\def \J{\text{\rm Jac}}
\def\jet#1#2{J^{#1}_{#2}}
\def\gp#1{G^{#1}}
\def\gpo{\gp {2k_0}_0}
\def\emmp {\scrF(M,p;M',p')}
\def\rk{\text{\rm rk\,}}
\def\Orb{\text{\rm Orb\,}}
\def\Exp{\text{\rm Exp\,}}
\def\Span{\text{\rm span\,}}
\def\d{\partial}
\def\D{\3J}
\def\pr{{\rm pr}}

\def \CZZ {\C \dbl Z,\zeta \dbr}
\def \D{\text{\rm Der}\,}
\def \Rk{\text{\rm Rk}\,}
\def \CR{\text{\rm CR}}
\def \ima{\text{\rm im}\,}
\def \I {\mathcal I}

\def \M {\mathcal M}

\newtheorem{Thm}{Theorem}[section]
\newtheorem{Cor}[Thm]{Corollary}
\newtheorem{Pro}[Thm]{Proposition}
\newtheorem{Lem}[Thm]{Lemma}

\theoremstyle{definition}\newtheorem{Def}[Thm]{Definition}

\theoremstyle{remark}
\newtheorem{Rem}[Thm]{Remark}
\newtheorem{Exa}[Thm]{Example}
\newtheorem{Exs}[Thm]{Examples}

\numberwithin{equation}{section}

\def\bl{\begin{Lem}}
\def\el{\end{Lem}}
\def\bp{\begin{Pro}}
\def\ep{\end{Pro}}
\def\bt{\begin{Thm}}
\def\et{\end{Thm}}
\def\bc{\begin{Cor}}
\def\ec{\end{Cor}}
\def\bd{\begin{Def}}
\def\ed{\end{Def}}
\def\be{\begin{Exa}}
\def\ee{\end{Exa}}
\def\bpf{\begin{proof}}
\def\epf{\end{proof}}
\def\ben{\begin{enumerate}}
\def\een{\end{enumerate}}

\newcommand{\dbar}{\bar\partial}
\newcommand{\genmat}{\lambda}
\newcommand{\polynorm}[1]{{|| #1 ||}}
\newcommand{\vnorm}[1]{\left\|  #1 \right\|}
\newcommand{\asspol}[1]{{\mathbf{#1}}}
\newcommand{\Cn}{\mathbb{C}^n}
\newcommand{\Cd}{\mathbb{C}^d}
\newcommand{\Cm}{\mathbb{C}^m}
\newcommand{\C}{\mathbb{C}}
\newcommand{\CN}{\mathbb{C}^N}
\newcommand{\CNp}{\mathbb{C}^{N^\prime}}
\newcommand{\Rd}{\mathbb{R}^d}
\newcommand{\Rn}{\mathbb{R}^n}
\newcommand{\RN}{\mathbb{R}^N}
\newcommand{\R}{\mathbb{R}}
\newcommand{\bR}{\mathbb{R}}
\newcommand{\N}{\mathbb{N}}
\newcommand{\dop}[1]{\frac{\partial}{\partial #1}}
\newcommand{\vardop}[3]{\frac{\partial^{|#3|} #1}{\partial {#2}^{#3}}}
\newcommand{\br}[1]{\langle#1 \rangle}
\newcommand{\infnorm}[1]{{\left\| #1 \right\|}_{\infty}}
\newcommand{\onenorm}[1]{{\left\| #1 \right\|}_{1}}
\newcommand{\deltanorm}[1]{{\left\| #1 \right\|}_{\Delta}}
\newcommand{\omeganorm}[1]{{\left\| #1 \right\|}_{\Omega}}
\newcommand{\nequiv}{{\equiv \!\!\!\!\!\!  / \,\,}}
\newcommand{\bk}{\mathbf{K}}
\newcommand{\p}{\prime}
\newcommand{\tV}{\mathcal{V}}
\newcommand{\poly}{\mathcal{P}}
\newcommand{\ring}{\mathcal{A}}
\newcommand{\ringk}{\ring_k}
\newcommand{\ringktwo}{\mathcal{B}_\mu}
\newcommand{\germs}{\mathcal{O}}
\newcommand{\On}{\germs_n}
\newcommand{\mcl}{\mathcal{C}}
\newcommand{\formals}{\mathcal{F}}
\newcommand{\Fn}{\formals_n}
\newcommand{\autM}{{\Aut (M,0)}}
\newcommand{\autMp}{{\Aut (M,p)}}
\newcommand{\holmaps}{\mathcal{H}}
\newcommand{\biholmaps}{\mathcal{B}}
\newcommand{\autmaps}{\mathcal{A}(\CN,0)}
\newcommand{\jetsp}[2]{ G_{#1}^{#2} }
\newcommand{\njetsp}[2]{J_{#1}^{#2} }
\newcommand{\jetm}[2]{ j_{#1}^{#2} }
\newcommand{\glnc}{\mathsf{GL_n}(\C)}
\newcommand{\glmc}{\mathsf{GL_m}(\C)}
\newcommand{\glc}{\mathsf{GL_{(m+1)n}}(\C)}
\newcommand{\glk}{\mathsf{GL_{k}}(\C)}
\newcommand{\smC}{\mathcal{C}^{\infty}}
\newcommand{\anC}{\mathcal{C}^{\omega}}
\newcommand{\kC}{\mathcal{C}^{k}}



\begin{abstract} We prove the following CR version of Artin's approximation theorem for holomorphic mappings between real-algebraic sets in complex space. Let $M\subset \C^N$ be a real-algebraic CR submanifold whose CR orbits are all of the same dimension.  Then for every point $p\in M$, for every real-algebraic subset $S'\subset \C^N\times\C^{N'}$ and every positive integer $\ell$, if $f\colon (\C^N,p)\to \C^{N'}$ is a germ of a holomorphic map such that ${\rm Graph}\, f \cap (M\times \C^{N'})\subset  S'$, then there exists a germ of a complex-algebraic map $f^\ell \colon (\C^N,p)\to \C^{N'}$ such that ${\rm Graph}\, f^\ell \cap (M\times \C^{N'})\subset  S'$ and that agrees with $f$ at $p$ up to order $\ell$. 

\vskip 0.5\baselineskip

\noindent{\sc R\'esum\'e.} Nous d\'emontrons la version CR suivante du th\'eor\`eme d'approximation d'Artin pour des applications holomorphes entre sous-ensembles alg\'ebriques r\'eels des espaces euclidiens complexes. Soit $M\subset \C^N$ une sous-vari\'et\'e CR alg\'ebrique r\'eelle dont les orbites CR sont toutes de m\^eme dimension. Pour tout point $p\in M$, pour tout sous-ensemble alg\'ebrique r\'eel  $S'\subset \C^N\times\C^{N'}$ et pour tout entier naturel $\ell$, si $f\colon (\C^N,p)\to \C^{N'}$ est un germe d'application holomorphe tel que ${\rm Graph}\, f \cap (M\times \C^{N'})\subset  S'$, alors il existe un germe d'application alg\'ebrique complexe  $f^\ell \colon (\C^N,p)\to \C^{N'}$ telle que  ${\rm Graph}\, f^\ell \cap (M\times \C^{N'})\subset  S'$ et dont le jet d'ordre  $\ell$ en $p$ coincide avec celui de $f$. 
\end{abstract}

\maketitle

\section{Introduction}\Label{int}

A well-known theorem in algebraic geometry going back to Artin \cite{A69} states that given any system of polynomial equations $P(x,y)=0$, $x=(x_1,\ldots,x_n)$, $y=(y_1,\ldots,y_m)$ over the field of real numbers (resp.\  complex numbers) and given any germ of an analytic solution $y(x)$ of the above system at a given point $p\in \R^n$ (resp.\ $\C^n$), there exists a sequence of germs at $p$ of {\em real-algebraic} (resp.\ {\em complex-algebraic}) solutions of the system that converge to the given solution in the Krull topology. 

In this paper, we provide a Cauchy-Riemann version of this approximation theorem. Given a real-algebraic CR submanifold $M\subset \C^N$, let us say that $M$ has the {\em Nash-Artin approximation property} if for every point $p\in M$, every real-algebraic subset $S'\subset \C^{N+N'}$ and every positive integer $\ell$, if $f\colon (\C^N,p)\to \C^{N'}$ is a germ of a  holomorphic map such that ${\rm Graph}\, f \cap (M\times \C^{N'})\subset  S'$ (as germs at $(p,f(p))$), then there exists a germ at $p$ of complex-algebraic map $f^\ell \colon (\C^N,p)\to \C^{N'}$ that agrees with $f$ at $p$ up to order $\ell$ with  ${\rm Graph}\, f^\ell \cap (M\times \C^{N'})\subset  S'$.  Observe that finding a sequence $(f^\ell)_\ell$ with the above properties consists of finding a sequence of germs of real-algebraic mappings from $M$ into $\C^{N'}\simeq \R^{2N'}$ whose graphs are contained in $S'\subset \C^{N+N'}\simeq \R^{2(N+N')}$ but which, furthermore, must satisfy the tangential Cauchy-Riemann equations (on $M$). Hence, deciding whether $M$ has the Nash-Artin approximation property is indeed equivalent to asking for a CR (or more generally PDE) version of Artin's theorem mentioned above.

We provide a positive solution to this approximation problem by making use of the CR geometry of the manifold $M$.  Let $T^{c}M\subset  TM$ denote the complex tangent bundle of $M$. Recall that for every $p\in M$, there exists a unique germ of a  real-algebraic CR submanifold $\mathcal{O}_p$ through $p$ with 
the property that every point $q\in \mathcal{O}_p$ can be reached from $p$ by following a piecewise differentiable curve  in $M$ whose tangent vectors are 
in $T^cM$  (see \cite{BER96}). We call this germ  {\em the CR orbit} of $M$ at $p$  and  say that $M$ is {\em minimal} at $p$ if this CR orbit is a neighbhorhood of $p$ in $M$. We also say  that $M$ is of {\em constant orbit dimension} if the CR orbits of $M$ have all the same dimension. Our main result is the following:

\begin{Thm}\label{t:nageneralbis} Let $M\subset \C^N$ be a real-algebraic CR submanifold of constant orbit dimension. Then $M$ has the Nash-Artin approximation property.
\end{Thm}

The study of algebraicity properties of holomorphic mappings sending real-algebraic sets into each other has been extensively studied in the past years (see e.g.\ \cite{W77, H94, HJ98, BER96, M98, CMS99, Z99, BRZ2, BMR, MMZsurvey, LM09}). Most of the mentioned results are concerned with the automatic algebraic extension of all holomorphic mappings between two given real-algebraic sets and hold under some geometric nondegeneracy conditions on these sets. The situation considered in Theorem~\ref{t:nageneralbis} goes beyond this setting and deals with the general case where $M$ is allowed to be mapped to an arbitrary real-algebraic set by a non-algebraic holomorphic map. In such a situation, the only known sufficient condition implying that a (connected) real-algebraic CR submanifold possesses the Nash-Artin approximation property is that of minimality and goes back to the work of Meylan, Zaitsev and the author \cite{MMZ3, MMZsurvey}. Theorem~\ref{t:nageneralbis} provides the same conclusion under a much weaker sufficient condition, that, in addition, is generically satisfied on every connected real-algebraic CR manifold (see e.g.\ \cite{BERbook}). Hence even the following immediate consequence of Theorem~\ref{t:nageneralbis} is new.


\begin{Cor}\label{c:subvariety}
For every connected real-algebraic CR submanifold $M\subset \C^N$, there exists a closed proper real-algebraic subvariety $\Sigma_M$ of $M$  such that $M\setminus \Sigma_M$ has the Nash-Artin approximation property.
\end{Cor}





In the case of real hypersurfaces, we can deduce from Theorem~\ref{t:nageneralbis} the following stronger result.

\begin{Thm}\label{t:naghyp}
Any $($smooth$)$ real-algebraic hypersurface of $\C^N$ has the Nash-Artin approximation property.
\end{Thm}

Indeed, any connected component of any real-algebraic (smooth) hypersurface of $\C^N$ is either Levi-flat, in which case it is of constant orbit dimension and Theorem~\ref{t:nageneralbis} applies, or is somewhere minimal, in which case it has the Nash-Artin approximation property in view of \cite[Theorem 3.5]{MMZsurvey}. 


One noteworthy and typical application of Theorem~\ref{t:nageneralbis} deals with the problem of deciding whether the holomorphic equivalence of two germs of real-algebraic CR manifolds of the same dimension implies their algebraic equivalence (see the works \cite{BRZ2, BMR, LM09}). Specializing Theorem~\ref{t:nageneralbis} to this situation, we have:

\begin{Cor}\label{c:alghol}
 Let $M,M'\subset \C^N$ be two real-algebraic CR submanifolds of constant orbit dimension. Then  for all points $p\in M$ and $p'\in M'$, the holomorphic equivalence of the germs $(M,p)$ and $(M',p')$ implies their algebraic equivalence.
\end{Cor}

As for Theorem~\ref{t:nageneralbis}, Corollary~\ref{c:alghol} was known only when $M$ and $M'$ are (connected and somewhere)  minimal and goes back in this case from the work of Baouendi, Rothschild and the author \cite{BMR}. For nowhere minimal real-algebraic CR manifolds, partial results towards Corollary~\ref{c:alghol} were previously established by Baouendi, Rothschild and Zaitsev \cite{BRZ2} and more recently by Lamel and the author in \cite{LM09}. In these two works, the conclusion given by Corollary~\ref{c:alghol} is obtained for all points $p$ in a certain Zariski open subset of $M$. The proofs of \cite{BRZ2, LM09} require to exclude from $M$ a thin set of points corresponding to the locus of the degeneracy set of a certain holomorphic foliation. Corollary~\ref{c:alghol} answers one of the main questions left open from \cite{BRZ2, LM09} which was to decide whether one could get rid off this locus.

Let us now discuss briefly the proof of Theorem~\ref{t:nageneralbis}. If $M$ is as in Theorem~\ref{t:nageneralbis} and $p\in M$, we assume that $p=0$ and we may view the germ of $M$ at $0$ as a (small) algebraic deformation of its CR orbits. Namely, there
exists an integer $c\in \{0,\ldots,N\}$ and a real-algebraic submersion $S \colon (M,0) \to (\R^c,0)$ such 
that $S^{-1} (S(q)) = \mathcal{O}_q $ for all $q\in M$ near $0$. The level sets of $S$, i.e.\  the CR orbits, therefore foliate $M$ near the origin by {\em minimal} real-algebraic CR submanifolds. We may hence identify the germ of $M$ at $0$ with an algebraic deformation $(M_t)$ for $t\in  \R^c$ close to $0$, where $M_t\subset \C^{N-c}$ is a germ at $0$ of a {\em minimal} real-algebraic CR submanifold (see e.g.\ Lemma~\ref{l:coordconstantorbit}). If $f\colon (\C^{N-c}_z\times \C^c_u,0)\to \C^{N'}$ is a germ of a holomorphic map such that ${\rm Graph}\, f \cap (M\times \C^{N'})\subset  S'$ for some real-algebraic subset $S'\subset \C^N \times \C^{N'}$, then for every $t\in \R^c$ sufficiently small, the holomorphic map $f_t\colon (\C^{N-c}_z,0)\ni z \mapsto f(z,t)\in  \C^{N'}$ satisfies ${\rm Graph}\, f_t \cap (M_t\times \C^{N'})\subset  S'_t$ where $S_t':=\{(z,z')\in \C^{N-c} \times \C^{N'}:(z,t,z')\in S'\}$. Since each submanifold $M_t$ is minimal, the conclusion of  Theorem~\ref{t:nageneralbis} boils down to providing a deformation version of the approximation theorem of \cite{MMZ3, MMZsurvey} associated to the deformation $(M_t)_{t\in \R^c}$, the analytic family of holomorphic maps $(f_t)_{t\in \R^c}$ and the family of real-algebraic subsets $(S_t')_{t\in \R^c}$.

To this end, we devote one part of the paper, namely Section~\ref{s:exam}, to build a system of polynomial equations with real-algebraic coefficients fulfilled by the restrictions to $\R^N$ and $\R^c$ of $f$ and $J(u):=(\partial^{\alpha}f(0,u);|\alpha|\leq k)$ for some suitable integer $k$. The system is constructed in such a way that it reduces our approximation problem to an application of a version of Artin's approximation theorem due to Popescu \cite{P} (see \cite{LM09} where the use of \cite{P} appeared for the first time in the subject). To construct the desired system, we introduce a suitable field of partially algebraic power series, denoted by $\4S^f$, depending on the mapping $f$ and its derivatives. An appropriate and careful study of the field extension generated by the field $\4S^f$ and the components of the mapping $f$ leads to the desired system of real-algebraic polynomial equations. At this point, we should mention that the above strategy depends on a key technical proposition, Proposition~\ref{p:keyprop}, which is an algebraic dependence result for certain power series which might be of independent interest. We devote Section~\ref{s:crratios} to the proof of this result : it requires to introduce algebraic power series rings defined over a field that is itself a field of fractions of partially algebraic power series and depending on  the fixed mapping $f$. During the proof of Proposition~\ref{p:keyprop}, we adapt several tools concerning ratios of power series developed in \cite{MMZ3} and also make use of Artin's approximation theorem over the above defined field of partially algebraic power series.  We should point out that, under the additional assumption that the source manifold $M$ in Theorem~\ref{t:nageneralbis} is connected and contains minimal points, most of the tools developed in this paper are unecessary and the proof can be simplified and reduced to that of \cite{MMZ3, MMZsurvey}.

The paper is organized as follows. In Section~\ref{s:preliminaries}, we introduce the notation and preliminary notions that will be used throughout the article. Section~\ref{s:crratios} is devoted to the proof of one key algebraicity property for certain ratios of power series that is used in Section~\ref{s:exam} while building the system of real-algebraic equations associated to a given mapping $f$. The proof of Theorem~\ref{t:nageneralbis} is completed in Section~\ref{s:final}.

\section{Preliminaries}\label{s:preliminaries}

\subsection{Algebraic power series and two approximation theorems for polynomial systems} 
Throughout the paper, given a field $\4K$, and indeterminates $x=(x_1,\ldots,x_r)$, $r\geq 1$, we denote by $\4K [[x]]$ the ring of formal power series with coefficients in the ring $\K$. For $T(x)\in \4K [[x]]$ and an integer $m$, we write $j^mT$ or  $(j^mT)(x)$ for the collection of all partial derivatives of $T$ up to order $m$. If the indeterminates $x$ split as $x=(x',x'')$, we write  $(j^mT)(0,x'')$ for the power series mapping $(j^mT)$ evaluated at $x'=0$. We will also be using the notation $j_{x'}^mT$ for the collection of all partial derivatives with respect to $x'$ up to order $m$.

We shall frequently make use of the ring of algebraic power series over $\4K$ : it is the subring $\6A_{\K}\{x\}$ of $\K [[x]]$   consisting of those $T(x)\in \K [[x]]$ for which there exist a positive integer $m$ and, for $j=0,\dots,m$, polynomials $P_j \in \4K[x]$  with $P_m \neq 0$ such that \[ \sum_{j=0}^m P_j (x) T(x)^j = 0. \]

Note that when $\4K=\C$, $\6A_{\C}\{x\}$ is the usual ring of (complex-)algebraic (or Nash) functions over $\C$, which is contained in the ring of convergent power series denoted in this paper by $\C \{x\}$. Given any convergent power series
$\eta=\eta (x)\in \C\{x\}$, we will also denote throughout the paper by $\bar{\eta}=\bar{\eta}(x)$ the convergent power series obtained from $\eta$ by taking complex conjugates of its coefficients.

In this work, we will make use of two approximation theorems for polynomial systems of equations with algebraic coefficients. The first one is due to Artin \cite{A69}.

\begin{Thm}\label{t:artin}$(${\rm Artin \cite{A69}}$)$ Let $\4K$ be a field, $P(x,y)=(P_1(x,y),\ldots,P_m(x,y))$ be $m$ polynomials in the ring $\6A_{\K}\{x\}[y]$ with $x=(x_1,\ldots,x_r)$, $y=(y_1,\ldots,y_q)$. Suppose that $Y(x)=(Y_1(x),\ldots,Y_q(x))\in (\4K [[x]])^q$ is a formal solution of the system
\begin{equation}\label{e:system}
 P(x,Y(x))=0.
\end{equation}
Then, for every integer $\ell$, there exists $Y^\ell(x)\in (\6A_{\K}\{x\})^q$ satisfying the system \eqref{e:system} such that $Y^\ell(x)$ agrees with $Y(x)$ up to order $\ell$.
\end{Thm}

The second one is a more precise version of Theorem~\ref{t:artin} and corresponds to Popescu's solution of the approximation problem on nested subrings.

\begin{Thm}\label{t:popescu}$(${\rm Popescu} \cite{P}$)$ Let $\4K$ be a field, $P(x,y)=(P_1(x,y),\ldots,P_m(x,y))$ be $m$ polynomials in the ring $\6A_{\K}\{x\}[y]$ with $x=(x_1,\ldots,x_r)$, $y=(y_1,\ldots,y_q)$. Suppose that $Y(x)=(Y_1(x),\ldots,Y_q(x))\in (\4K [[x]])^q$ is a formal solution of the system
\begin{equation}\label{e:systembis}
 P(x,Y(x))=0.
\end{equation}
Suppose furthermore that for every $j=1,\ldots,q$, there exists $s_j\in \{1,\ldots,r\}$ such that $Y_j(x)\in \4K [[x_1,\ldots,x_{s_j}]]$. Then, for every integer $\ell$, there exists $Y^\ell(x)\in (\K\{x\})^q$ satisfying the system \eqref{e:systembis}, such that $Y^\ell(x)$ agrees with $Y(x)$ up to order $\ell$ and  each $Y^\ell_j(x)\in \6A_{\4K} \{x_1,\ldots,x_{s_j}\}$, $j=1,\ldots,q$.
\end{Thm}

Observe that by a noetherianity argument, Theorems~\ref{t:artin} and \ref{t:popescu} also hold for a polynomial system with infinitely many equations.

\subsection{$\4K$-algebraicity of ratios of formal power series}
In Section~\ref{s:crratios}, we will have to study certain types of ratios of formal power series defined over a ground field that is not the usual field of complex numbers, but, rather a field extension over $\C$. In order to prove a key property of these ratios, we shall use and modify some concepts about ratios of formal power series introduced in \cite{MMZ3} for the case of $\C[[x]]$ to treat ratios of power series over $\4K[[x]]$ where $\C \hookrightarrow \4K$ is a field extension. 

Let  $\C \hookrightarrow \4K$ be a field extension. Given two power series $N(x),D(x)\in \4K[[x]]$, $x=(x_1,\ldots,x_r)$, we write $(N:D)$ for a pair of formal power series where  we both allow $N$ and $D$ to be zero. When $D\not \equiv 0$, we can think of $(N:D)$ as being the usual ratio $N/D$. In the following, we say that two ratios $(N_1:D_1)$ and $(N_2:D_2)$ of formal power series in $\4K [[x]]$ are equivalent if $N_1D_2-N_2D_1=0$. As in \cite{MMZ3}, it will be useful to introduce the following notion.
 
\begin{Def}  Given two ratios $(N_1:D_1)$ and $(N_2:D_2)$ of formal power series in $\4K [[x]]$, a complex linear subspace $E\subset \C^r$ and a nonnegative integer $q$, we say that $(N_1:D_1)$ and $(N_2:D_2)$ are $q$-similar along $E$ if $(j^q(N_1D_2-N_2D_1))|_E=0$.
\end{Def}

We have the following (weak) transitivity property of the notion of $q$-similarity, whose simple proof is completely analogous to that of \cite[Lemmas 3.2 and 3.3]{MMZ3} and is left to the reader.

\begin{Lem}\label{l:transitivity}
Let $(N_1:D_1)$, $(N_2:D_2)$ and $(N_3:D_3)$ be three ratios of formal power series in $\4K [[x]]$. Suppose that we have a splitting of the indeterminates of the form $x=(u_1,u_2,u_3)$ and let $E$ be the linear subspace of $\C^r$ given by $E=\{u_1=0,\, u_2=0\}$. Suppose that there exist an integer $\ell$ and integers $q_1,q_2\geq \ell$ such that 
$$(j^\ell(N_2,D_2))|_E\not = 0,\quad (j_{u_1}^{q_1}j^{q_2}_{u_2}(N_1D_2-N_2D_1))|_E = 0,$$
$$(j_{u_1}^{q_1}j^{q_2}_{u_2}(N_3D_2-N_2D_3))|_E = 0.$$
Then $(j_{u_1}^{q_1-\ell}j^{q_2-\ell}_{u_2}(N_3D_1-N_1D_3))|_E=0$.
\end{Lem}

We now define the notion of $\4K$-algebraicity for (formal) power series along a certain complex subspace that will be extremely useful in Section~\ref{s:crratios}.

\begin{Def}\label{d:kalgebraic}  Given a ratio $(N:D)$ of formal power series in $\4K [[x]]$ and a complex linear subspace $E\subset \C^r$, we say that $(N:D)$ is $\4K$-algebraic along $E$ if there exists an integer $\ell$ and for every integer $q$, formal power series $N_q,D_q\in  \6A_{\K}\{x\}$ such that $(N:D)$ and $(N_q:D_q)$ are $q$-similar along $E$ and such that $(j^\ell(N_q,D_q))|_E\not \equiv 0$.
\end{Def}

\begin{Rem}\label{r:bof}
 Using Lemma~\ref{l:transitivity}, it is easy to see that if two ratios of formal power series in $\4K[[x]]$ are equivalent and if one of them is nontrivial and $\4K$-algebraic along a complex subspace $E\subset \C^r$ then the other ratio has the same property.
\end{Rem}

\subsection{Generic real-algebraic submanifolds and local CR orbits}\label{ss:iterated} 
Let $M\subset \C^N$ be a real-algebraic CR submanifold with $N\geq 2$ and $T^{0,1}M$ its CR bundle. For every point $p\in M$, we denote by ${\mathcal G}_M(p)$ the Lie algebra evaluated at $p$ generated by the sections of $T^{0,1}M$ and its conjugate $T^{1,0}M$. By a theorem of Nagano (see e.g.\ \cite{BERbook, BCH08}), for every point $p\in M$, there is a well-defined unique germ at $p$ of a real-analytic submanifold $V_p$ satisfying $\C T_qV_p={\mathcal G}_M(q)$ for all $q\in V_p$. This unique submanifold is necessarily CR and is called the {\em CR orbit} of $M$ at $p$. In fact, by \cite[Corollary 2.2.5]{BER96}, this CR orbit is even a {\em real-algebraic} CR submanifold contained in $M$. It is not difficult to see that if $M$ is connected, the dimension of the local CR orbits is constant (and of maximal dimension) except possibly on a proper  real-algebraic subvariety $\Sigma_M$ of $M$ (see e.g.\ \cite{BRZ1, BRZ2}). If ${\rm dim}\, V_p={\rm dim}\, M$ for some point $p\in M$, we say that $M$ is of  {\em minimal} (or also of {\em finite type}) at $p$.

Assume in what follows that $M$ is a (connected) real-algebraic generic submanifold of $\C^N$ of CR dimension $n$ and codimension $d$. If $p$ is a point in $M\setminus \Sigma_M$, where $\Sigma_M$ is defined as above, then one may choose holomorphic (algebraic) coordinates such that $M$ is described near $p$ through the following well-known lemma.

\begin{Lem}\label{l:coordconstantorbit} $($\cite[Proposition 3.4]{BRZ1}\, {\rm and}\, \cite[Lemma 3.4.1]{BER96}$)$ 
	Let $M\subset \C^N$ be a connected generic real-algebraic submanifold through a point $p\in M$ whose CR orbit at $p$ is of maximal dimension and let $c\in \{0,\ldots, d\}$ be the codimension of this CR orbit in $M$. 
	Then there exists normal algebraic coordinates $Z=(z,\eta)\in \C^n \times \C^d$, $\eta=(w,u)\in  \C^{d-c}\times \C^c$, 
	such that $M$ is given near the origin by an equation of the form
\begin{equation}\label{e:conormal}
 \eta=(w,u)=\Theta (z,\bar z,\bar \eta):=(Q(z,\bar z, \bar w,\bar u),\bar u), 
\end{equation}
where $Q$ is a $\C^{d-c}$-valued complex-algebraic map near $0\in \C^{n+N}$.
Furthermore, there exist 
neighborhoods $U,V$ of the origin in $\R^c$ and $\C^{N-c}$ respectively such that for every $u\in U$, the real-algebraic 
submanifold given by 
\begin{equation}\label{e:Mu}
M_u:=\{(z,w)\in V: w=Q(z,\bar z, \bar w,u)\}
\end{equation}
is generic in $\C^{N-c}$ and minimal at $0$.
\end{Lem}

As a real-analytic submanifold, $M$ can be complexified and gives rise to its so-called  {\em complexification}, which 
we 
denote by $\6M$. This complexification $\6M$ is the germ at $0$ of the complex-algebraic submanifold of $\C^{2N}$ given by 
\begin{equation}\label{e:complexification}
\{(Z,\zeta)\in (\C^N\times \C^N,0):\sigma =\bar \Theta (\chi,z,\eta)\},
\end{equation}
where $Z=(z,\eta)\in \C^n \times \C^d$ and $\zeta=(\chi,\sigma)\in \C^n \times \C^d$. Recall also that the above choice of normal coordinates imply that the following two identities hold
\begin{equation}\label{e:normalreal}
 \Theta (z,0,\sigma)=\Theta (0,\chi,\sigma)=\sigma,\quad \quad \Theta (z,\chi,\bar \Theta (\chi,z,\eta))=\eta.
\end{equation}

In what follows, we pick a point $p\in M\setminus \Sigma_M$ and choose normal coordinates for $M$ at $p$ as given by Lemma~\ref{l:coordconstantorbit}. 
In Section~\ref{s:crratios}, we will need the {\em iterated Segre mappings} attached to such a germ of submanifold 
(see e.g.\ \cite{BRZ1}). These mappings are defined in the following way. For any nonnegative integer $j$, 
 we  denote by $t^j$ a variable lying in $\C^{n}$ and also introduce the variable 
$t^{[j]}:=(t^1,\ldots,t^j)\in \C^{nj}$. We also use the notation $s$ for a variable lying in the euclidean space $\C^{d-c}$. Then we first set 
$v_0(s,u):=(0,s,u)$ for $(s,u)\in \C^d$ sufficiently close to $0$ and define 
the  map $v_j\colon (\C^{nj}\times \C^{d-c}\times \C^c,0)\to \C^N$ for $j\geq 1$ inductively  as follows:
\begin{equation}\label{e:iteratedmaps}
v_j(t^{[j]},s,u):=(t^j, U_j(t^{[j]},s,u)),\quad {\rm where}\quad U_j(t^{[j]},s,u):=\Theta (t^j, \bar v_{j-1}(t^{[j-1]},s,u)).
\end{equation}
From the construction we clearly see that each iterated Segre mapping $v_j$ defines an algebraic map in a neighbhorhood of $0$ in $\C^{nj+d}$.  In fact, for every point $(s,u)\in \C^d$ sufficiently close to the origin, the map $v_j(\cdot,s,u)$ parametrizes the usual Segre set of order $j$ attached to the point $(0,s,u)$ (see e.g.\ \cite{BER96, BERbook}). Notice also that, thanks to \eqref{e:normalreal}, one has the following useful  identities
\begin{equation}\label{e:identities}
v_j(0,s,u)=(0,s,u),\quad v_{j+2}(t^{[j+2]},s,u)|_{t^{j+2}=t^{j}}=v_{j}(t^{[j]},s,u),\quad j\geq 0, 
\end{equation}
and that for every $j\geq 0$, the germ at $0$ of the holomorphic map $(v_{j},\bar v_{j-1})$ takes its values in $\6M$, where $\6M$ is the complexification of $M$ as defined by \eqref{e:complexification}.

\section{An algebraicity property for certain ratios of power series on real-algebraic generic submanifolds}\label{s:crratios}
Throughout this section, we assume that $M$ is a germ of a (connected) real-algebraic generic submanifold through the origin in $\C^N$ with $N\geq 2$ and that $f\colon (\C^N,0)\to \C^{N'}$ is a germ of a holomorphic mapping. We also assume that the CR orbit of $M$ at the origin is of maximal dimension and that normal coordinates $Z=(z,w,u)$ for $(M,0)$ have been chosen (and fixed) as in Lemma~\ref{l:coordconstantorbit}. In such a setting, we are going to prove a crucial algebraicity property for certain ratios of (convergent) power series constructed from the mapping $f$ and whose restriction on $M$ is CR. 

\subsection{Statement of the algebraicity property} Given any integer $k$, we denote by $J_0^k(\C^N,\C^{N'})$ the jet space of order $k$ at the origin of holomorphic maps from $\C^N$ to $\C^{N'}$. Throughout the paper, $\Lambda^k$, $\Gamma^k$ and $\Upsilon^k$ will denote coordinates in $J_0^k(\C^N,\C^{N'})$ and we write $\Lambda^k=(\Lambda_\alpha^k)_{|\alpha|\leq k}$ where $\Lambda_\alpha^k \in \C^{N'}$ for $\alpha \in \N^N$ (and analogously for $\Gamma^k$ and $\Upsilon^k$).

In order to state the algebraic criterion given by Proposition~\ref{p:keyprop} below, we need to define the following ring of power series that depend on the above fixed choice of normal coordinates $Z=(z,w,u)$.

\begin{Def}\label{d:bf} Let $f$ and $Z=(z,w,u)$ be as above. Let $\4B^f$ be the subring of $\C \{u\}[z,w]$ consisting of those power series $T(z,w,u)$ for which there exists an integer $k$ and $S\in \C[z,w,u,\Lambda^k,\Gamma^k]$ such that $T(z,w,u)=S(z,w,u,(j^kf)(0,u), (j^k\bar f)(0,u))$. The field of fractions  of $\4B^f$ will be denoted by $\4S^f$.
\end{Def}

In what follows, we say that a holomorphic vector $X$ defined  near $0\in \C^N_Z \times \C^N_\zeta$ is a $(0,1)$ vector field if it annihilates the natural projection $\C^N \times \C^N \ni (Z,\zeta)\mapsto Z\in \C^N$. The goal of this section is to prove the following:

\begin{Pro}\label{p:keyprop} Let $M$, $f$ and $Z=(z,w,u)$ be as above. Let $\phi_1,\phi_2 \in \C \{Z,\zeta\}$ with $\phi_2|_{\6M} \not \equiv 0$, where $\6M$ is the complexification of $M$ as given by \eqref{e:complexification}. Assume that there exist an integer $k$, and power series $\Phi_1,\Phi_2 \in \6A_{\C}\{Z,\zeta\}[\Lambda^k,\Gamma^k,\Upsilon^k]$ such that
\begin{equation}\label{e:phi}
\phi_j (Z,\zeta)= \Phi_j (Z,\zeta,(j^kf)(0,u), (j^k\bar f)(0,u),(j^k \bar f)(\zeta))),\ j=1,2.
\end{equation}
Assume furthermore that the ratio $\phi_1/\phi_2$ is annihilated by any $(0,1)$ holomorphic vector field tangent to $\6M$ near $0$. Then $\phi_1/\phi_2$ is algebraic over the field ${\4S^f}$ $($as introduced in Definition~{\rm \ref{d:bf}}$)$.
\end{Pro}

Other (different) types of CR ratios of (formal) power series on real-analytic generic submanifolds have been studied in \cite{MMZ3}. Though the technique used in \cite{MMZ3} deals only minimal generic submanifolds, we will  use part of this technique  to prove Proposition~\ref{p:keyprop}.

\subsection{Proof of Proposition~\ref{p:keyprop}}
We assume in what follows that we are in the setting of Proposition~\ref{p:keyprop}. In order to prove Proposition~\ref{p:keyprop}, we need to work in appropriate rings of formal power series. All these formal power series will have the same ground field that is defined as follows.

\begin{Def}\label{d:kf} Let $\4D_u^f$ be the subring of $\C \{u\}$ consisting of those (convergent) power series $T(u)$ for which there exists an integer $k$ and $S\in \6A_{\C}\{u\}[\Lambda^k,\Gamma^k]$ such that $T(u)=S(u,(j^kf)(0,u), (j^k\bar f)(0,u))$. The field of fractions of the ring $\4D_u^f$ will be denoted by $\4K_u^f$. 
\end{Def}

The first step in the proof of Proposition~\ref{p:keyprop} relies on the following lemma.

\begin{Lem}\label{l:stepone}
 In the setting of Proposition~{\rm \ref{p:keyprop}}, there exists power series $\Delta_1 (Z),\Delta_2 (Z)\in \C \{Z\}\cap \K_u^f[[z,w]]$ such that $\Delta_2 \not \equiv 0$ and such that the ratios $(\phi_1:\phi_2)$ and $(\Delta_1:\Delta_2)$ are equivalent. In addition, the ratio $(\Delta_1:\Delta_2)$ $($viewed as a ratio of power series in $\K_u^f[[z,w]])$ is $\4K_u^f$-algebraic along the $(n+c)$-dimensional complex subspace $\{w=0\}$ $($as explained in Definition~{\rm \ref{d:kalgebraic}}$)$.
\end{Lem}
\begin{proof} By assumption, we know that \begin{equation}\label{e:derive}
\frac{\partial}{\partial \chi}\left(\frac{\phi_1 (Z,\chi,\bar \Theta (\chi, Z))}{\phi_2 (Z,\chi,\bar \Theta (\chi, Z))}\right)=0.
\end{equation}
Since $\phi_2|_{\6M} \not =0$, we can choose a multiindex $\beta \in \N^n$ of minimal length such that $$\left(\partial_\chi^{|\beta|}\left(\phi_2 (Z,\chi,\bar \Theta (\chi, Z))\right)\right)\Big|_{\chi=0}\not \equiv 0.$$
Setting $$\Delta_1 (Z):=\left(\partial_\chi^{|\beta|}\left(\phi_1 (Z,\chi,\bar \Theta (\chi, Z))\right)\right)\Big|_{\chi=0},\quad \Delta_2 (Z):=\left(\partial_\chi^{|\beta|}\left(\phi_2 (Z,\chi,\bar \Theta (\chi, Z))\right)\right)\Big|_{\chi=0},$$ 
it follows from \eqref{e:derive} that $(\Delta_1:\Delta_2)$ is a ratio of convergent power series that is equivalent to $(\phi_1:\phi_2)$. In addition, it follows from the chain rule and from \eqref{e:phi} that there exists power series $\widetilde \Phi_j \in \6A_{\C}\{Z\}[\Lambda^k,\Gamma^k,\Upsilon^{k+|\beta|}]$, $j=1,2$, such that
\begin{equation}\label{e:phi2}
\Delta_j (Z)= \widetilde\Phi_j (Z,(j^kf)(0,u), (j^k\bar f)(0,u),(j^{k+|\beta|} \bar f)(0,w,u))),\ j=1,2.
\end{equation}
Then \eqref{e:phi2} readily implies that both power series $\Delta_1,\Delta_2$ belong to  $\K_u^f[[z,w]]$. In fact, one may even  notice that for every integer $q$, if $\Delta_{1,q}$ and $\Delta_{2,q}$ denote the truncated power series
$$\Delta_{1,q}(z,w,u)=\sum_{\nu \in \N^{d-c}\atop |\nu|\leq q}\frac{\partial^\nu \Delta_1}{\partial w^\nu}(z,0,u)\, \frac{w^\nu}{\nu !},\quad \Delta_{2,q}(z,w,u)=\sum_{\nu \in \N^{d-c}\atop |\nu|\leq q}\frac{\partial^\nu \Delta_2}{\partial w^\nu}(z,0,u)\, \frac{w^\nu}{\nu !},$$
then $\Delta_{1,q},\Delta_{2,q} \in  \6A_{\4K_u^f}\{z,w\}$, and $(\Delta_{1,q}:\Delta_{2,q})$ and $(\Delta_1:\Delta_2)$ are $q$-similar along the subspace $L:=\{w=0\}$. Furthermore, since there exists an integer $n_0$ such that $(j^{n_0}(\Delta_1,\Delta_2))\Big|_{L}~\not \equiv 0$, it follows that for all $q\geq n_0$, we have $(j^{n_0}(\Delta_{1,q},\Delta_{2,q}))\Big|_{L}\not \equiv 0$, which shows that the ratio $(\Delta_1,\Delta_2)$  is $\4K_u^f$-algebraic along the complex subspace $L$. The proof of Lemma~\ref{l:stepone} is complete.
\end{proof}

For the second step of the proof of Proposition~\ref{p:keyprop}, we shall use the iterated Segre maps introduced in section~\ref{ss:iterated}. For any integer $j$, recall that $v_j$ is the mapping given by \eqref{e:iteratedmaps}. Note that since each map $v_j$ is algebraic and since the power series $\Delta_1=\Delta_1 (z,w,u)$ and $\Delta_2=\Delta_2 (z,w,u)$ given by Lemma~\ref{l:stepone} belong to the ring $\K_u^f[[z,w]]$, it follows from the first identity of \eqref{e:identities} that the power series $\Delta_1 \circ v_j=(\Delta_1 \circ v_j)(t^{[j]},s,u)$ and $\Delta_2 \circ v_j=(\Delta_2 \circ v_j)(t^{[j]},s,u)$ both belong to the ring  $\4K_u^f[[t^{[j]},s]]$. In what follows, when writing any iterated map $v_j$, we shall sometimes omit to write the variables for sake of brevity.

\begin{Lem}\label{l:iterationproc}
 In the setting of Proposition~{\rm \ref{p:keyprop}}, let $j$ be a positive integer and assume that the ratio $(\Delta_1 \circ v_j:\Delta_2 \circ v_j)$ of power series in $\4K_u^f[[t^{[j]},s]]$ is $\4K_u^f$-algebraic along the complex subspace $E_j:=\{s=0\}\subset \C^{nj+d-c}$. Then the ratio $(\Delta_1 \circ v_{j+2}:\Delta_2 \circ v_{j+2})$ of power series in $\4K_u^f[[t^{[j+2]},s]]$ is $\4K_u^f$-algebraic along the complex subspace $E_{j+2}:=\{s=0\}\subset \C^{n(j+2)+d-c}$.
\end{Lem}

In order to prove Lemma~\ref{l:iterationproc}, we need the following preliminary result.
\begin{Lem}\label{l:new}
Under the assumptions of Lemma~{\rm \ref{l:iterationproc}}, the ratio $(\Delta_1 \circ v_{j+2}:\Delta_2 \circ v_{j+2})$ is $\4K_u^f$-algebraic along the complex subspace $\widetilde E_{j+2}:=\{s=0,\ t^{j+2}=t^j\}\subset \C^{n(j+2)+d-c}$.
\end{Lem}

\begin{proof}[Proof of Lemma~{\rm \ref{l:new}}]
By assumption, there exist an integer $\ell_0$ and for every integer $q$, (formal) power series $N_q,D_q\in \6A_{\4K_u^f}\{t^{[j]},s\}$ such that 
\begin{equation}\label{e:startup}
 (j^{\ell_0}(N_q,D_q))\big|_{E_j}\not \equiv 0,\quad (j^{q}((\Delta_2 \circ v_j)N_q-(\Delta_1 \circ v_j)D_q))\big|_{E_j}=0.
\end{equation}
Consider the algebraic map $I\colon (\C^{nj+d},0)\ni (t^{[j]},s,u)\mapsto \C^{nj+d}$ given by
\begin{equation}\label{e:inversion}
I(t^{[j]},s,u):=
\begin{cases}
\begin{aligned}
(t^{[j]},U_j(t^j,t^{j-1},\ldots,t^1,s,u)) \quad {\rm for}\ j\ {\rm odd},\\
(t^{[j]},\overline{U}_j(t^j,t^{j-1},\ldots,t^1,s,u)) \quad {\rm for}\  j\  {\rm even},
\end{aligned}
\end{cases}
\end{equation}
where $U_j$ is the mapping given by \eqref{e:iteratedmaps}. The reader can check that, by using the second identity of \eqref{e:normalreal}, the map $I$ satisfies the following properties
\begin{equation}\label{e:properties}
 v_j(I(t^{[j]},s,u))=(t^j,s,u),\quad I(t^{[j-1]},v_j(t^{[j]},s,u))=(t^{[j]},s,u).
\end{equation}
Consider the algebraic map $J(t^{[j+2]},s,u):=I(t^{[j-1]},v_{j+2}(t^{[j+2]},s,u))$ and, for every integer $q$, define $\widetilde N_q:=N_q \circ J$ and $\widetilde D_q:=D_q\circ J$. Since $J(0,0,u)=(0,0,u)$ and since $J\in \6A_{\C}\{t^{[j+2]},s,u\}$, the power series $\widetilde N_q$ and $\widetilde D_q$ belong to the ring $\4K_u^f[[t^{[j+2]},s]]$. In fact, noticing that the last $\C^c$-valued component of $J(t^{[j+2]},s,u)$ is equal to $u$ and using the fact that $N_q,D_q\in \6A_{\4K_u^f}\{t^{[j]},s\}$ and the second identity of \eqref{e:properties}, we see that $\widetilde N_q,\widetilde D_q\in \6A_{\4K_u^f}\{t^{[j+2]},s\}$. Furthermore, it follows from \eqref{e:properties} that $v_{j+2}=v_j\circ J$ and therefore we have
\begin{equation}\label{e:onemore}
(\Delta_2 \circ v_{j+2})\widetilde N_q-(\Delta_1 \circ v_{j+2})\widetilde D_q=\left((\Delta_2 \circ v_j)N_q-(\Delta_1 \circ v_j)D_q\right)\circ J.
\end{equation}
Noticing furthermore that \eqref{e:properties} and \eqref{e:identities} imply that $J(\widetilde E_{j+2})\subset E_j$, it follows from \eqref{e:startup} and \eqref{e:onemore} that for every integer $q$, the ratios $(\Delta_1 \circ v_{j+2}:\Delta_2 \circ v_{j+2})$ and $(\widetilde N_q:\widetilde D_q)$ are  $q$-similar along $\widetilde E_{j+2}$. In addition, the second identity in \eqref{e:properties} implies $N_q=\widetilde N_q\big|_{t^{j+2}=t^j}$ and $D_q=\widetilde D_q\big|_{t^{j+2}=t^j}$ which shows, in view of the first identity in \eqref{e:startup} that $(j^{\ell_0}(\widetilde N_q,\widetilde D_q))\big|_{\widetilde E_{j+2}}\not \equiv 0$. This proves Lemma~\ref{l:new}.
\end{proof}

\begin{proof}[Proof of Lemma~{\rm \ref{l:iterationproc}}]
We first note that in view of \eqref{e:phi}, there exists power series $\Phi_{1,j+2}$ and $\Phi_{2,j+2}$ both in the ring $\6A_{\C}\{t^{[j+2]},s,u\}[\Lambda^k,\Gamma^k,\Upsilon^k]$ such that for $r=1,2$ one has

\begin{equation}\label{e:phi3}
\begin{aligned}
\phi_{r,j+2}(t^{[j+2]},s,u)&:=\phi_r (v_{j+2},\bar v_{j+1})\\
&=\Phi_{r,j+2} (t^{[j+2]},s,u,(j^kf)(0,u), (j^k\bar f)(0,u),(j^k \bar f)\circ \bar v_{j+1})).
\end{aligned}
\end{equation}


Note in addition that the (convergent) power series $\phi_{r,j+2}$, $r=1,2$, both belong to the ring $\4K_u^f[[t^{[j+2]},s]]$. Using the fact that $\phi_2 \big |_{\6M}\not \equiv 0$ and that the mapping $(v_{j+2},\bar v_{j+1})$ is clearly submersive from $\C^{n(j+2)+d}$ into $\C^N$, we get that there exists an integer $\ell_1$ such that \begin{equation}\label{e:nondegenerate}
(j^{\ell_1}(\phi_{1,j+2},\phi_{2,j+2}))\big|_{\widetilde E_{j+2}}\not \equiv 0.
\end{equation}
Furthermore, since the ratios $(\phi_1:\phi_2)$ and $(\Delta_1:\Delta_2)$ are equivalent by Lemma~\ref{l:stepone}, it follows that the ratios $(\phi_{1,j+2}:\phi_{2,j+2})$ and $(\Delta_1 \circ v_{j+2}:\Delta_2 \circ v_{j+2})$ are also equivalent. Hence from Lemma~\ref{l:new} and Remark~\ref{r:bof}, we get that the ratio $(\phi_{1,j+2}:\phi_{2,j+2})$ is $\4K_u^f$-algebraic along $\widetilde E_{j+2}$. We are now going to follow the strategy of the proof of \cite[Lemma~3.8]{MMZ3}. There exists an integer $\ell_2$  and for every integer $q$, power series $A_q, B_q\in \6A_{\4K_u^f}\{t^{j+2},s\}$ such that the ratios $(A_q: B_q)$ and  $(\phi_{1,j+2}:\phi_{2,j+2})$ are $q$-similar along $\widetilde E_{j+2}$ with
\begin{equation}\label{e:acta}
(j^{\ell_2}(A_q,B_q))\big|_{\widetilde E_{j+2}}\not \equiv 0.
\end{equation}
Furthermore, we may assume without loss of generality that $\ell_2\geq \ell_1$. For the rest of the proof, we fix an arbitrary integer $m\geq \ell_2$ and set $Y(t^{[j+1]},s,u):=(j^k \bar f)\circ \bar v_{j+1}$.  Then for every integer $q\geq m+\ell_2$, we have
\begin{equation}\label{e:barbe}
(j_{t^{j+2}}^{q-m}\, j_{s}^m(A_q\phi_{2,j+2}-\phi_{1,j+2}B_q))\big|_{\widetilde E_{j+2}}=0.
\end{equation}
Since $A_q$ and $B_q$ are both in the ring $\6A_{\4K_u^f}\{t^{[j+2]},s\}$ and in view of \eqref{e:phi3}, we may rewrite \eqref{e:barbe} as follows
\begin{equation}\label{e:postes}
\6S_{q,m} (t^{[j+1]},u,(j_s^mY)(t^{[j+1]},0,u))=0,\ \forall q\geq m+\ell_2,
\end{equation}
where $\6S_{q,m}$ is a polynomial in its last argument with coefficients in the ring $\6A_{\4K_u^f}\{t^{[j+1]}\}$. By construction, the power series mapping $(j_s^mY)(t^{[j+1]},0,u)$ belongs to the ring $\4K_u^f[[t^{j+1}]]$ and is a formal solution of the polynomial system of equations with coefficients in $\6A_{\4K_u^f}\{t^{[j+1]}\}$ provided by \eqref{e:postes}. By Theorem~\ref{t:artin}, for every integer $r$, we may find a power series mapping $Y^r(t^{[j+1]},u)\in \6A_{\4K_u^f}\{t^{[j+1]}\}$ such that 
\begin{equation}\label{e:zainab}
(j^r_{t^{[j+1]}}Y^r)\big|_{t^{[j+1]}=0}=\left(j^r_{t^{[j+1]}}\left((j_s^mY)\big|_{s=0}\right)\right)\big|_{t^{[j+1]}=0}
\end{equation}
 and satisfying 
\begin{equation}\label{e:hope}
 \6S_{q,m} (t^{[j+1]},u,Y^r(t^{[j+1]},u))=0,\ \forall q\geq m+\ell_2.
\end{equation}
Choose $T^r\in \6A_{\K_u^f}\{t^{[j+1]}\}[s]$ such that $(j^m_s T^r)\big|_{s=0}=Y^r$ and set for $e=1,2$
\begin{equation}\label{e:phir}
 \phi^r_{e,j+2}(t^{[j+2]},s,u):= \Phi_{e,j+2} (t^{[j+2]},s,u,(j^kf)(0,u), (j^k\bar f)(0,u),T^r(t^{[j+1]},s,u))).
\end{equation}
In view of \eqref{e:nondegenerate}, \eqref{e:zainab} and the fact that $m\geq \ell_2\geq \ell_1$, we may choose $r$ such that \begin{equation}\label{e:choice}
(j^{\ell_1}(\phi^r_{1,j+2},\phi^r_{2,j+2}))\big|_{\widetilde E_{j+2}}\not \equiv 0,
\end{equation}
and keep such a choice of $r$ for the remainder of this proof. From the construction of the power series mapping $T^r$, it follows that $\phi_{1,j+2}^r,\phi_{2,j+2}^r \in \6A_{\K_u^f}\{t^{[j+2]}\}$ and that the following identity holds
\begin{equation}\label{e:barbebis}
(j_{t^{j+2}}^{q-m}\, j_{s}^m(A_q\phi^r_{2,j+2}-\phi^r_{1,j+2}B_q))\big|_{\widetilde E_{j+2}}=0,\ \forall q\geq m+\ell_2.
\end{equation}
Lemma~\ref{l:transitivity} together with \eqref{e:barbe}, \eqref{e:barbebis} and \eqref{e:acta} implies
\begin{equation}\label{e:good}
 (j_{t^{j+2}}^{q-m-\ell_2}\, j_{s}^{m-\ell_2}(\phi_{1,j+2}\, \phi^r_{2,j+2}-\phi^r_{1,j+2}\,\phi_{2,j+2}))\big|_{\widetilde E_{j+2}}=0,\ \forall q\geq m+\ell_2.
\end{equation}
Hence, \eqref{e:good} obviously implies
\begin{equation}\label{e:goodbis}
 (j_{s}^{m-\ell_2}(\phi_{1,j+2}\, \phi^r_{2,j+2}-\phi^r_{1,j+2}\,\phi_{2,j+2}))\big|_{E_{j+2}}=0,
\end{equation}
which shows that the ratios $(\phi_{1,j+2}:\phi_{2,j+2})$ and $(\phi^r_{1,j+2}:\phi^r_{2,j+2})$ are $(m-\ell_2)$-similar along $E_{j+2}$. Since \eqref{e:choice} immediately implies that $(j^{\ell_1}(\phi^r_{1,j+2},\phi^r_{2,j+2}))\big|_{E_{j+2}}\not \equiv 0$, we obtain that the ratio $(\phi_{1,j+2}:\phi_{2,j+2})$ is $\4K_u^f$-algebraic along $E_{j+2}$. Since the ratio $(\Delta_1 \circ v_{j+2}:\Delta_2 \circ v_{j+2})$ is equivalent to $(\phi_{1,j+2}:\phi_{2,j+2})$, the conclusion of Lemma~\ref{l:iterationproc} follows from Remark~\ref{r:bof}.
\end{proof}

\begin{proof}[Completion of the proof of Proposition~{\rm \ref{p:keyprop}}]
In view of Lemma~\ref{l:stepone}, we have to show that $\Delta_1/\Delta_2$ is algebraic over the field $\4S^f$. We start by noticing that Lemma~\ref{l:stepone} provides the fact the ratio $(\Delta_1 \circ v_1:\Delta_2 \circ v_1)$ of power series in $\4K_u^f [[t^1,s]]$ is $\4K_u^f$-algebraic along the subspace $E_1:=\{s=0\}$. Applying Lemma~\ref{l:iterationproc}, we get that for every odd integer, the ratio  $(\Delta_1 \circ v_j:\Delta_2 \circ v_j)$ (of power series in the ring $\4K_u^f[[t^{[j]},s]]$) is $\4K_u^f$-algebraic along the complex subspace $E_j:=\{s=0\}\subset \C^{nj+d-c}$. Choose $j=2d+3$ for the remainder of this proof, where we recall that $d$ is the codimension of $M$ in $\C^N$. There exists an integer $\ell_3$ and for every integer $q$, power series $N_q,D_q\in \6A_{\4K_u^f}\{t^{[2d+3]},s\}$ such that the ratios $(N_q:D_q)$ and $(\Theta \circ v_{2d+3}:\Delta \circ v_{2d+3})$ are $q$-similar  along the subspace $\{s=0\}$ with 
\begin{equation}\label{e:done}
 (j^{\ell_3}(N_q,D_q))\big|_{s=0}\not \equiv 0.
\end{equation}
Take $q=\ell_3$ and choose $\beta_0\in  \N^{(2d+3)n+d-c}$ with $|\beta_0|\leq \ell_3$ of minimal length such that $$(\partial^{\beta_0} N_{\ell_3},\partial^{\beta_0} D_{\ell_3})\big|_{s=0}\not \equiv 0.$$ 
Then we have
\begin{equation}\label{e:lastid}
\left((\Delta_2 \circ v_{2j+3})\, \partial^{\beta_0}N_{\ell_3}-(\Delta_1 \circ v_{2j+3})\, \partial^{\beta_0}D_{\ell_3}\right)\big|_{s=0}=0.
\end{equation}
Since the manifold $M_0$ given by Lemma~\ref{l:coordconstantorbit} is minimal at the origin, the minimality criterion given in \cite{BER96, BERbook} shows that the mapping $\C^{(2d+3)n+c}\ni (t^{[2d+3}],u)\mapsto v_{2d+3}(t^{[2d+3]},0,u)\in \C^N$ is of generic rank $N$. Hence \eqref{e:lastid} implies that $(\partial^{\beta_0}D_{\ell_3})\big|_{s=0}\not \equiv 0$ and that the following identity holds in the quotient field of $\4K_u^f[[t^{2d+3}]]$  
\begin{equation}\label{e:quotient}
\frac{(\Delta_1 \circ v_{2j+3})(t^{[2d+3]},0,u)}{(\Delta_2 \circ v_{2j+3})(t^{[2d+3]},0,u)}=\frac{(\partial^{\beta_0}N_{\ell_3})(t^{[2d+3]},0,u)}{(\partial^{\beta_0}D_{\ell_3})(t^{[2d+3]},0,u)}.
\end{equation}
Since $N_{\ell_3},D_{\ell_3}\in \6A_{\4K_u^f}\{t^{[2d+3]},s\}$, it follows from \eqref{e:quotient} that the germ of the meromorphic function given by the left hand side of \eqref{e:quotient} is also algebraic over the quotient field of the ring $\4K_u^f[t^{[2d+3]}]$. Therefore there exist two positive integers $a,b$ and, for $\mu=0,\ldots,a$, polynomials $P_\mu=P_\mu (u,t^{[2d+3]},\Lambda^b,\Gamma^b)\in \6A_{\C}\{u\}[t^{[2d+3]},\Lambda^b,\Gamma^b]$ such that $$P_a (u,t^{[2d+3]},(j^bf)(0,u),(j^b\bar f)(0,u))\not \equiv 0$$ and such that the following identity holds
\begin{equation}\label{e:soso}
 \sum_{\mu=0}^aP_\mu \left(u,t^{[2d+3]},(j^bf)(0,u),(j^b\bar f)(0,u)\right) \left(\frac{(\Delta_1\circ v_{2j+3})(t^{[2d+3]},0,u)}{(\Delta_2 \circ v_{2j+3})(t^{[2d+3]},0,u)} \right)^\mu=0.
\end{equation}
In fact, the above mentioned minimality criterion  even provides points arbitrarily close to $0$ in $\C^{2d+3}$ such that the (algebraic) map $t^{[2d+3]}\mapsto v_{2d+3}(t^{[2d+3]},0,0)$ sends those points to the origin and has rank $N-c$ there. Pick a point $T^0$ with this property such that the left hand side of \eqref{e:quotient} is meromorphic in an open neighborhood $V\times W$ of the origin in $\C^{(2d+3)n}\times \C^c$ with $T^0 \in V$. From the rank theorem, we may find an algebraic map $\lambda \colon (\C^N,0)\to (\C^{(2d+3)n},T^0)$ such that $v_{2d+3}(\lambda (z,w,u),0,u)=(z,w,u)$. Composing \eqref{e:soso} with the algebraic mapping $\lambda$, we get
\begin{equation}\label{e:sososo}
 \sum_{\mu=0}^aP_\mu \left(u,\lambda (z,w,u),(j^bf)(0,u),(j^b\bar f)(0,u)\right) \left(\frac{\Delta_1 (z,w,u)}{\Delta_2 (z,w,u)} \right)^\mu=0,
\end{equation}
for $(z,w,u)\in \C^N$ sufficiently close to the origin. It is easy to see that one may choose the mapping $\lambda$ so that $P_a \left(u,\lambda (z,w,u),(j^bf)(0,u),(j^b\bar f)(0,u)\right)\not \equiv 0$, which proves that $\Delta_1/\Delta_2$ is algebraic over the quotient field of $\6A_{\4K_u^f}\{z,w\}$. Since this latter field is algebraic over the quotient field of $\4K_u^f[z,w]$, which is itself algebraic over the field $\4S^f$, we see that the proof of Proposition~\ref{p:keyprop} is complete.
\end{proof}

\subsection{An application of the algebraicity property} We shall now provide an application of Proposition~\ref{p:keyprop} in the spirit  of \cite[Proposition 4.3]{MMZ3}. This result will be one of the key point of the proof of Theorem~\ref{t:nageneralbis}.

\begin{Pro}\label{p:buffet} Let $f\colon (\C^N,0)\to \C^{N'}$ be a germ of a holomorphic mapping and $M$ be a connected real-algebraic generic submanifold through the origin. Assume that the CR orbit of $M$ at $0$ is of maximal dimension and choose normal coordinates $Z=(z,w,u)$ for $M$ near $0$ as in Lemma~{\rm \ref{l:coordconstantorbit}}. Let $\6M$ be the complexification of $M$ as defined in {\rm \eqref{e:complexification}} and assume the mapping $f$ splits as follows $f=(g,h)\in \C^{e}\times \C^{N'-e}$ for some integer $e\in \{1,\ldots,N\}$. Assume also that there exists  an integer $k_0$ and a polynomial $R\in \6A_{\C}\{Z,\zeta\}[\Lambda^{k_0},\Gamma^{k_0}][\xi',\varpi']$ where $(\xi',\varpi')\in \C^{e} \times \C^{e}$ such that:
\begin{enumerate}
 \item[(i)] $R(Z,\zeta,(j^{k_0}f)(0,u),(j^{k_0}\bar f)(0,u),\xi',\varpi') \not \equiv 0$  $(Z,\zeta,\xi',\varpi')\in \6M \times \C^{e}\times \C^{e}$ near $0$,
\item[(ii)] $R(Z,\zeta,(j^{k_0}f)(0,u),(j^{k_0}\bar f)(0,u),g(Z),\bar g (\zeta)) =0$ for $(Z,\zeta)\in \6M$ near $0$.
\end{enumerate}
Then  the components of the mapping $g$ are algebraically dependent over the field $\4S^f$ $($as introduced in Definition~{\rm \ref{d:bf}}$)$.
\end{Pro}

\begin{proof}Let $\6E$ be the set of all polynomials $B\in \C \{Z,\zeta\}[\xi']$ such that there exist an integer $m$ and $D\in \6A_{\C}\{Z,\zeta\}[\Lambda^{m},\Gamma^{m},\Upsilon^m][\xi'] $ such that $$B(Z,\zeta,\xi')=D(Z,\zeta,(j^{m}f)(0,u),(j^{m}\bar f)(0,u),(j^m \bar f)(\zeta),\xi'),$$ with $B(Z,\zeta,\xi')\not \equiv 0$ for $(Z,\zeta,\xi')\in \6M \times \C^{N'}$ near $0$ and such that $B(Z,\zeta,g(Z))=0$ for $(Z,\zeta)\in \6M$ near $0$.

We claim that $\6E$ is not empty. Indeed, let $R$ as by the assumption of the lemma. There are two cases to consider. On one hand, if $$R(Z,\zeta,(j^{k_0}f)(0,u),(j^{k_0}\bar f)(0,u),\xi',\bar g (\zeta)) \not \equiv 0,\quad (Z,\zeta,\xi')\in \6M \times \C^{e}\  {\rm near}\ 0,$$ then we are clearly done. On the other hand, if for $(Z,\zeta,\xi')\in \6M \times \C^{e}$ near $0$ $$R(Z,\zeta,(j^{k_0}f)(0,u),(j^{k_0}\bar f)(0,u),\xi',\bar g (\zeta))=0,$$ then we also have $$\bar R(\zeta,Z,(j^{k_0}\bar f)(0,u),(j^{k_0} f)(0,u),\varpi', g (Z))=0$$ for $(Z,\zeta,\varpi')\in \6M \times \C^{e}$. Writing 
$$\bar R(\zeta,Z,\Gamma^{k_0},\Lambda^{k_0},\varpi',\xi')=\sum_{\alpha \in \N^{N'}}C_\alpha (\zeta,Z,\Gamma^{k_0},\Lambda^{k_0},\xi')(\varpi')^\alpha,$$
where the above sum is finite and each $C_\alpha \in \6A_{\C}\{Z,\zeta\}[\Gamma^{k_0},\Lambda^{k_0}][\xi']$,  assumption (i) implies that there exists $\beta \in \N^{N'}$ such that $C_\beta (\zeta,Z,(j^{k_0}\bar f)(0,u),(j^{k_0} f)(0,u),\xi')\not \equiv 0$ for $(Z,\zeta,\xi')\in \6M \times \C^{e}$ near $0$. Since  $C_\beta (\zeta,Z,(j^{k_0}\bar f)(0,u),(j^{k_0} f)(0,u),g(Z))= 0$ for $(Z,\zeta)\in \6M$ near $0$, the claim is proved.

Since $\6E$ is not empty, we may choose $B_0\in \6E$ with the minimum number $r_0$ of monomials in $\xi'$. If $r_0=1$, we see that one component of the mapping $g$ necessarily vanishes identically, which proves the proposition in this case. Otherwise,  let $m\in \4Z_+$ such that we may write $B_0(Z,\zeta,\xi')=D_0(Z,\zeta,(j^{m}f)(0,u),(j^{m}\bar f)(0,u),(j^m \bar f)(\zeta),\xi')$ for some $D_0\in \6A_{\C}\{Z,\zeta\}[\Lambda^{m},\Gamma^{m},\Upsilon^m][\xi']$. We expand $D_0$ in monomials in $\xi'$ in the following way
$$D_0(Z,\zeta,\Lambda^m,\Gamma^m,\Upsilon^m,\xi')=\sum_{j=1}^{r_0}D_{0,j}(Z,\zeta,\Lambda^m,\Gamma^m,\Upsilon^m)\, E_j(\xi').$$
We also set for $j=1\ldots,r_0$, \begin{equation}\label{e:Tj}
T_j(Z,\zeta):=D_{0,j}(Z,\zeta,(j^{m}f)(0,u),(j^{m}\bar f)(0,u),(j^m \bar f)(\zeta)),
\end{equation}
 and note that from the definition of $r_0$, it follows that $T_j(Z,\zeta)\not \equiv 0$ for $(Z,\zeta)\in \6M$ near $0$ and for every $j=1,\ldots,r_0$. Let $\6L_1,\ldots,\6L_n$ be a local basis of $(0,1)$ holomorphic  vector fields tangent to $\6M$ near $0$, where we recall that $n$ is the CR dimension of $M$. Since $M$ is real-algebraic, we may assume, without loss of generality, that these vector fields have (complex-)algebraic coefficients. Since   \begin{equation}\label{e:nom}
B_0(Z,\zeta,g(Z))=\sum_{j=1}^{r_0}T_j(Z,\zeta)\, E_j(g(Z))=0,\  (Z,\zeta)\in \6M \ {\rm near}\ 0,
\end{equation}
applying each vector $\6L_\nu$ to this identity yields,
\begin{equation}\label{e:lead}
 \sum_{j=1}^{r_0-1}\6L_\nu \left(\frac{T_j(Z,\zeta)}{T_{r_0}(Z,\zeta)}\right)E_j(g(Z))= 0,
\end{equation}
for $\nu=1,\ldots,n$ and $(Z,\zeta)\in \6M$ near $0$. For every $\nu \in \{1,\ldots,n\}$ and $j=1\ldots,r_0-1$, we note that the ratio of convergent power series $\6L_\nu \left(\frac{T_j(Z,\zeta)}{T_{r_0}(Z,\zeta)}\right)$ may be written in the form  $T_{j,\nu}(Z,\zeta)/(T_{r_0}(Z,\zeta))^2$ where
$$T_{j,\nu}(Z,\zeta)=G_{j,\nu}(Z,\zeta,(j^{m+1}f)(0,u),(j^{m+1}\bar f)(0,u),(j^{m+1} \bar f)(\zeta))$$ for some $G_{j,\nu}\in \6A_{\C}\{Z,\zeta\}[\Lambda^{m+1},\Gamma^{m+1},\Upsilon^{m+1}]$. Hence, it follows from \eqref{e:lead} and the minimality of $r_0$ that  $\6L_\nu \left(\frac{T_j(Z,\zeta)}{T_{r_0}(Z,\zeta)}\right)\equiv 0$ for all $\nu=1,\ldots,r_0-1$. We may therefore apply Proposition~\ref{p:keyprop} to conclude that each ratio $T_j/T_{r_0}$ is algebraic over the field $\4S^f$ (as given in Definition~\ref{d:bf}). In other words, we may say that in the quotient field of $\C\{Z,\zeta\}$, the subfield $\4F$ generated by $\4S^f$ and all the ratios $T_j/T_{r_0}$, $j=1\ldots,r_0-1$, is an algebraic extension of $\4S^f$. Since \eqref{e:nom} tells us that the components of the mapping $g$ are algebraically dependent over $\4F$, it follows that they are also algebraically dependent over  $\4S^f$. The proof is therefore complete.
\end{proof}

\section{A polynomial system of equations with real-algebraic coefficients associated to a given holomorphic map}\label{s:exam}
In this section, we assume that we are in the same situation as in the previous section. Namely, we assume that $M$ is a germ of a (connected) real-algebraic generic submanifold through the origin in $\C^N$ with $N\geq 2$ and that $f\colon (\C^N,0)\to \C^{N'}$ is a germ of a holomorphic mapping. We also assume that the CR orbit of $M$ at $0$ is of maximal dimension and that normal coordinates $Z=(z,w,u)\in \C^n \times \C^{d-c}\times \C^c$ for $(M,0)$ have been chosen (and fixed) as in Lemma~\ref{l:coordconstantorbit}.

Our goal in this section is to construct a system of real-algebraic equations over $\R^N$ associated to the mapping $f$ such that $f|_{\R^N}$ and $(j^kf)|_{\{0\}\times \R^c}$ are solutions of such a system for a suitable integer $k$. The constructed system will possess further additional properties (see Proposition~\ref{p:summarize} for the precise statement) that will be crucial in order to prove Theorem~\ref{t:nageneralbis}.

\subsection{Field extension associated to a given holomorphic map} We refer the reader to \cite{ZS58} for the basic notions of field theory used in the remainder of this paper and keep using the notation introduced in previous sections.

Let $\4S^f(f)$ denote the field generated by $\4S^f$ and all the components of the mapping $f$. Then denoting $\4M \{Z\}$ the field of fractions of $\C\{Z\}$, we have the following field extensions $$\4S^f\hookrightarrow \4S^f(f) \hookrightarrow \4M \{Z\}.$$

Let $e=e(f)\in \{0,\ldots,N'\}$ be te the transcendence degree of the field extension $\4S^f\hookrightarrow \4S^f(f)$. We may therefore choose $e$ components of the mapping $f$, denoted by $g$ in the rest of this paper, such that $g$ forms a transcendence basis of $\4S^f(f)$ over $\4S^f$. 
\bigskip

{\em We shall assume in the rest of Section \ref{s:exam} that $e<N'$.}
\bigskip

If we write $f=(g,h)\in \C^e_{\xi'} \times \C^{N'-e}_{\omega'}$, then the components of $g$ are all algebraically independent over $\4S^f$ and any other component of $h$ is algebraically dependent over $\4S^f(g)$. Writing $h=(h_1,\ldots,h_{N'-e})$ and $\omega'=(\omega_1',\ldots,\omega_{N'-e}')$, we may find an integer $k$ and, for every $j\in \{1,\ldots,N'-e\}$, a polynomial $\6P_j\in \C[Z,\Lambda^k,\Gamma^k,\xi'][\omega_j']$ such that for all $Z=(z,w,u)\in \C^N$ sufficiently close to the origin
\begin{equation}\label{e:out}
 \6P_j(Z,(j^kf)(0,u),(j^k\bar f)(0,u),g(Z),h_j(Z))= 0.
\end{equation}
Furthermore, writing $\6P_j(Z,\Lambda^k,\Gamma^k,\xi',\omega_j')=\sum_{\nu =0}^{m_j}\6P_{j,\nu}(Z,\Lambda^k,\Gamma^k,\xi')(\omega_j')^\nu$, we may assume that 
\begin{equation}\label{e:assure}
 \6P_{j,m_j}\left(Z,(j^kf)(0,u),(j^k\bar f)(0,u),g(Z)\right)\not \equiv 0,\quad j=1,\ldots,N'-e.
\end{equation}
We now come to a useful lemma that is a more refined version of \cite[Lemma~6.2]{MMZ3} (or of \cite[Lemma~4.3]{LM09})  but whose proof is completely analogous. For sake of completeness, we provide the proof of this statement. In what follows, we write $Z'=(\xi',\omega')\in \C^{e}\times \C^{N'-e}$ and will denote by $\varpi'$ another variable lying in $\C^e$ and $\theta'=(\theta_1',\ldots,\theta_{N'-e}') \in \C^{N'-e}$. We also denote $\ell_k={\rm card}\, \{\alpha \in \N^N:|\alpha|\leq k\}$.

\begin{Lem}\label{l:miami} In the above setting, for every real-valued polynomial $\rho(Z,\bar Z,Z',\bar Z')$ in $\C^N \times \C^{N'}$, there exists a nontrivial polynomial $\6Q^{\rho}\in \C [Z,\zeta,\Lambda^k,\Gamma^k,\xi',\varpi'][X]$ such that the following holds: for all germs of holomorphic maps $H\colon (\C^N,0)\to \C^{N'-e}$,  $G\colon (\C^N,0)\to \C^{e}$,  $T \colon (\C^c,0)\to \C^{\ell_k}$ satisfying 
$$\6P_j(Z,T(u),\bar T(u),G(Z),H_j(Z))\equiv 0,\ j=1,\ldots,N'-e,$$
and
$$ \6P_{j,m_j}\left(Z,T(u),\bar T(u),G(Z)\right)\not \equiv 0,$$
then 
$$\6Q^{\rho} (Z,\zeta,T(u),\bar T(u),G(Z),\bar G(\zeta),\rho (H(Z),G(Z),\bar H (\zeta),\bar G(\zeta)))\equiv 0,$$
for all $(Z,\zeta)\in \C^{2N}$ near $0$. Furthermore, we may write 
$$\6Q^{\rho}(Z,\zeta,\Lambda^k,\Gamma^k,\xi',\varpi',X)=\sum_{\nu =0}^\delta \6Q^{\rho}_\nu (Z,\zeta,\Lambda^k,\Gamma^k,\xi',\varpi') X^\nu,$$
with $\delta$ independent of $\rho$ and with $\6Q^{\rho}_\delta (Z,\zeta,T(u),\bar T(u),G(Z),\bar G(\zeta))\not \equiv 0$ for $(Z,\zeta)\in \6M$ near the origin.
\end{Lem}

\begin{proof} For $j=1,\ldots,N'-e$, consider 
$$\6P_j(Z,\Lambda^k,\Gamma^k,\xi',\omega_j')=\sum_{\nu =0}^{m_j}\6P_{j,\nu}(Z,\Lambda^k,\Gamma^k,\xi')(\omega_j')^\nu,$$
$$ \6V_j(\zeta,\Lambda^k,\Gamma^k,\varpi',\theta_j'):=\sum_{\nu =0}^{m_j}\bar{\6P}_{j,\nu}(\zeta,\Gamma^k,\Lambda^k,\varpi')(\theta_j')^\nu.$$
For every $j\in \{1,\ldots,N'-e\}$, let $A_j$ (resp.\ $B_j$) be the complex-algebraic variety in $\C^{N+2\ell_k+e}$ given by the zero set of the polynomial $\6P_{j,m_j}$ (resp.\ $\bar{\6P}_{j,m_j}$) and set $E:=\cup_{j=1}^{N'-e}A_j\cup \cup_{j=1}^{N'-e}B_j$. For $(Z,\Lambda^k,\Gamma^k,\xi')\in \C^{N+2\ell_k+e} \setminus E$ (resp.\ $(\zeta,\Gamma^k,\Lambda^k,\varpi')\in \C^{N+2\ell_k+e} \setminus E$) and for every $j\in \{1,\ldots,N'-e\}$, denote by
 $\sigma^{(1)}_j(Z,\Lambda^k,\Gamma^k,\xi'),\ldots,\sigma^{(m_j)}_j(Z,\Lambda^k,\Gamma^k,\xi')$ (resp.\ $\varsigma^{(1)}_j(\zeta,\Gamma^k,\Lambda^k,\varpi'),\ldots,\varsigma^{(m_j)}_j(\zeta,\Gamma^k,\Lambda^k,\varpi')$) the $m_j$ roots of the polynomial $\6P_j$ (resp.\ $\6V_j$).

Next for $(Z,\Lambda^k,\Gamma^k,\xi')$ and $(\zeta,\Gamma^k,\Lambda^k,\varpi')$ as above, consider the following polynomial $W$ in $X$
\begin{multline}
W(Z,\zeta,\Lambda^k,\Gamma^k,\xi',\varpi',X):=\\
\prod_{i_1=1}^{m_1}\ldots \prod_{i_{N'-e}=1}^{m_{N'-e}}\prod_{n_1=1}^{m_1}\ldots \prod_{n_{N'-e}=1}^{m_{N'-e}}\left( X-\rho (Z,\zeta,\xi',\sigma_1^{(i_1)},\ldots,\sigma_{N'-e}^{(i_{N'-e})},\varpi',\varsigma_1^{(n_1)},\ldots,\varsigma_{N'-e}^{(n_{N'-e})})\right),
\end{multline}
 where, for $r=1,\ldots,N'-e$, we have written $\sigma_{r}^{(i_{r})}$ for $\sigma_{r}^{(i_{r})}(Z,\Lambda^k,\Gamma^k,\xi')$ and $\varsigma_r^{(n_r)}$ for $\varsigma_r^{(n_r)}(\zeta,\Gamma^k,\Lambda^k,\varpi')$.
It follows from Newton's theorem on symmetric polynomials that $W$ can be rewritten in the following form
\begin{equation}
W(Z,\zeta,\Lambda^k,\Gamma^k,\xi',\varpi',X)=
X^\delta +\sum_{\gamma<\delta}W_\gamma(Z,\zeta,\Lambda^k,\Gamma^k,\xi',\varpi')X^\gamma,
\end{equation}
where $\delta$ is independent of $\rho$ and where for each $\gamma=0,\ldots,\delta-1$, 
\begin{multline}W_\gamma (Z,\zeta,\Lambda^k,\Gamma^k,\xi',\varpi'):=\\D_\gamma \left(Z,\zeta,\xi',\varpi',\left( \left(\frac{\6P_{j,\nu_j}(Z,\Lambda^k,\Gamma^k,\xi')}{\6P_{j,m_j}(Z,\Lambda^k,\Gamma^k,\xi')}\right)_{ \nu_j }\right)_{j},\left( \left(\frac{\bar{\6P}_{j,\nu_j}(\zeta,\Gamma^k,\Lambda^k,\varpi')}{\bar{\6P}_{j,m_j}(\zeta,\Gamma^k,\Lambda^k,\varpi')}\right)_{ \nu_j}\right)_{j}\right),
\end{multline}
for some polynomial $D_\gamma$ of its arguments (depending only on $\rho$) and where $1\leq j\leq N'-e$ and $1\leq \nu_j \leq m_j$. Setting $$K(Z,\zeta, \Lambda^k,\Gamma^k,\xi',\varpi'):=\prod_{j=1}^{N'-e}\6P_{j,m_j}(Z,\Lambda^k,\Gamma^k,\xi')\cdot \bar{\6P}_{j,m_j}(\zeta,\Gamma^k,\Lambda^k,\varpi'),$$
we see that for a suitable integer $r_\rho$, $\6Q^{\rho}:=K^{r_\rho}W \in \C [Z,\zeta,\Lambda^k,\Gamma^k,\xi',\varpi'][X]$. We leave it to the reader to check that the construction of the  polynomial $\6Q^{\rho}$ provides all the desired properties of the lemma. The proof is therefore complete.
\end{proof}

\subsection{Construction of the polynomial system attached to the mapping $f$} 
For every real polynomial $\rho(Z,\bar Z, Z',\bar Z')$ defined over $\C^N \times \C^{N'}$, let $\6Q^{\rho}$ be the polynomial given by Lemma~\ref{l:miami}. As in Lemma~\ref{l:miami}, we write
\begin{equation}\label{e:neweq}
\6Q^{\rho}(Z,\zeta,\Lambda^k,\Gamma^k,\xi',\varpi',X)=\sum_{\nu =0}^{\delta} \6Q^{\rho}_\nu (Z,\zeta,\Lambda^k,\Gamma^k,\xi',\varpi') X^\nu,
\end{equation}
and define 
\begin{equation}\label{e:ictp}
p_{\rho}:={\rm inf}\, \{\nu \in \{0,\ldots,\delta\} :\6Q^{\rho}_\nu(Z,\zeta,(j^kf)(0,u),(j^k\bar f)(0,u),\xi',\varpi')\big|_{\6M \times \C^{2e}}\not \equiv 0,\ {\rm near}\ \, 0.\}.
\end{equation}
It follows from \eqref{e:out} and \eqref{e:assure} and Lemma~\ref{l:miami} that each $p_{\rho}$ is well defined.

We need the following lemma.

\begin{Lem}\label{l:systT} In the above setting, there exists a finite collection of polynomials $\6S_1,\ldots,\6S_b\in \6A_{\C}\{u\}[\Lambda^k,\Gamma^k]$ such that for every $T\in (\C\{u\})^{\ell_k}$, $T$ satisfies the system of equations
\begin{equation}\label{e:uir}
\6Q^{\rho}_\nu(Z,\zeta,T(u),\bar T(u),\xi',\varpi')\big|_{\6M \times \C^{2(N'-e)}}\equiv 0,\ {\rm near}\ \, 0,
\end{equation}
for all $ \nu \in \{ 0,\ldots,p_{\rho}-1\}$ and every real-valued polynomial $\rho$ defined over $\C^N \times \C^{N'}$
 if and only if 
$$\6S_q(u,T(u),\overline{T(u)})\equiv 0, \ {\rm near}\, \, 0\in \R^c,\, q=1,\ldots,b.$$
\end{Lem}

\begin{proof}
For every real-valued polynomial $\rho$ defined over $\C^N \times \C^{N'}$ and every $0\leq \nu \leq p_\rho -1$, we write
$$\6Q^{\rho}_\nu(Z,\zeta,\Lambda^k,\Gamma^k,\xi',\varpi')=\sum_{\alpha,\beta \in \N^{e}}\6Q^{\rho}_{\nu,\alpha,\beta}(Z,\zeta,\Lambda^k,\Gamma^k)(\xi')^\alpha (\varpi')^\beta.$$
Then a convergent power series mapping $T$ as in Lemma~\ref{l:systT} satisfies \eqref{e:uir} if and only if for every  $\alpha,\beta \in \N^{N'-e}$, every $\nu \in \{0,\ldots, p_\rho-1\}$, we have
\begin{equation}\label{e:firstpartsystem}
 \6Q^\rho_{\nu,\alpha,\beta}(Z,\zeta,T(u),\bar T (u))\equiv 0,\ (Z,\zeta)\in \6M,\ \, {\rm close}\, \,{\rm to}\, \, 0.
\end{equation}
Using Lemma~\ref{l:coordconstantorbit}, we may choose  a parametrization $(\R^{2N-d-c}\times \R^c,0)\ni (y,u)\mapsto \varphi (y,u)\in \C^N$ of $M$ near $0$ such that $\varphi$ is real-algebraic and such that for each fixed $u\in \R^c$ sufficiently close to the origin, the mapping $(\R^{2N-d-c},0)\ni y\mapsto \varphi (y,u)$ parametrizes (the germ at $0$ of) the manifold $M_u$ as defined in \eqref{e:Mu}. Hence for all $u\in \C^c$ sufficiently close to the origin and for every $\gamma \in \N^{2N-d-c}$ and every $\alpha,\beta,\nu$, we may set
\begin{equation}\label{e:ines}
 \6S^\rho_{\nu,\alpha,\beta,\gamma}(u,\Lambda^k,\Gamma^k):=\left( \frac{\partial^\gamma}{\partial y^\gamma}\left(\6Q^\rho_{\nu,\alpha,\beta}(\varphi (y,u),\bar{\varphi} (y,u),\Lambda^k,\Gamma^k)\right)\right)\Big|_{y=0}.
\end{equation}
From our construction, each $\6S^\rho_{\nu,\alpha,\beta,\gamma}\in \6A_{\C}\{u\}[\Lambda^k,\Gamma^k]$. Furthermore a convergent power series mapping $T(u)$ as in Lemma~\ref{l:systT} satisfies \eqref{e:firstpartsystem} if and only if  for every $\gamma \in \N^{2N-d-c}$ and every $\alpha,\beta,\nu$ as above 
$$\6S^\rho_{\nu,\alpha,\beta,\gamma}(u,T(u),\bar T (u))\equiv 0$$ for all $u\in \R^c$ sufficiently close to the origin. By a noetherian argument, we may extract a finite collection of polynomials in the family $(\6S^\rho_{\nu,\alpha,\beta,\gamma})$ providing the conclusion of the lemma. The proof is therefore complete.
\end{proof}

We are now in a position to prove the main technical result of this section.

\begin{Pro}\label{p:summarize} Let $M$ and $f$ be above. Let $e$ be the transcendence degree of the field extension $\4S^f\hookrightarrow \4S^f(f)$ and assume that $e<N'$. For $j\in \{1,\ldots,N'-e\}$ and $q\in \{1,\ldots,b\}$, let $\6P_j\in \C[Z,\Lambda^k,\Gamma^k,\xi'][\omega_j']$ and $\6S_q \in \6A_{\C}\{u\}[\Lambda^k,\Gamma^k]$ given by \eqref{e:out} and Lemma~{\rm \ref{l:systT}} respectively. For $x\in \R^{n+d-c}$, $u\in \R^c$, consider the following $($complex-valued$)$ polynomial system 
\begin{equation}\label{e:mainsystem}
 \6P_j(x,u,\Lambda^k,\overline{\Lambda}^k, \xi',\omega_j')=0,\quad \6S_q(u,\Lambda^k,\overline{\Lambda}^k)=0,\ j\in \{1,\ldots,N'-e\},\ q\in \{1,\ldots,b\}.
\end{equation}
Then $\Lambda^k=(j^kf)(0,u)$, $\xi'=g(x,u)$, $\omega'=(\omega_1',\ldots,\omega_{N'-e}')=h(x,u)$ is a complex-valued real-analytic solution of \eqref{e:mainsystem}. Furthermore, the system \eqref{e:mainsystem} has the following property: for every real-algebraic set $\Sigma'\subset \C^N \times \C^{N'}$ and every sequence $\Lambda^k=T^\ell(u)$, $\xi'=g^\ell(x,u)$, $\omega'=h^\ell(x,u)$ of germs at $0\in \R^N$ of real-analytic mappings converging as $\ell \to \infty$ in the Krull topology to $(j^kf)(0,u)$, $g(x,u)$ and $h(x,u)$ respectively and satisfying \eqref{e:mainsystem}, if ${\rm Graph}\, f \cap (M\times \C^{N'})\subset \Sigma'$ $($as germs at $(0,f(0)))$, then for $\ell$ large enough, we also have ${\rm Graph}\, f^\ell \cap (M\times \C^{N'})\subset \Sigma'$ $($as germs at $(0,f(0)))$ where $f^\ell\colon (\C^N,0)\to \C^{N'}$ is the holomorphic map obtained by complexifying the real-analytic mapping $(g^\ell,h^\ell)\colon (\R^N,0)\to \C^{N'}$.
\end{Pro}

\begin{proof} The first part of the conclusion of the proposition follows immediately the construction of the system \eqref{e:mainsystem} (e.g.\ Lemma~\ref{l:systT} and \eqref{e:ictp} and \eqref{e:out}). Let us the prove the second part of the desired conclusion. To this end, let $\Sigma'$ be a fixed real-algebraic subset of $\C^N \times \C^{N'}$ (that we may assume to be, without loss of generality, different from $\C^N \times \C^{N'}$) and assume that ${\rm Graph}\, f \cap (M\times \C^{N'})\subset \Sigma'$ $($as germs at $(0,f(0)))$. Let also $(T^\ell(u))_\ell$, $(g^\ell(x,u))_\ell$, $(h^\ell(x,u))_\ell$ be a sequence of germs at $0\in \R^N$ of real-analytic mappings converging as $\ell \to \infty$ in the Krull topology to $(j^kf)(0,u)$, $g(x,u)$ and $h(x,u)$ respectively, and, satisfying for every integer $\ell$ the system \eqref{e:mainsystem}. We may assume that $\Sigma'\not = \C^N \times \C^{N'}$ since otherwise the conclusion of the proposition is obvious.

We prove the proposition by contradiction. Assume therefore that there exists a subsequence $(f^{\eta_r})_r$ such that  ${\rm Graph}\, f^{\eta_r} \cap (M\times \C^{N'})\not \subset \Sigma'$ for every integer $r$. Choose a finite number of nontrivial real-valued polynomials $\rho_1,\ldots,\rho_m$ defined over $\C^N \times \C^{N'}$ such that $\Sigma'=\cap_{j=1}^m \Sigma_j'$ where $\Sigma_j'$ is the zero set of the polynomial $\rho_j$. Then, by the pigeonhole principle, we may assume that there exists a subsequence $(f^{\tilde \eta_r})_r$ of  $(f^{\eta_r})_r$ such that  ${\rm Graph}\, f^{\tilde \eta_r} \cap (M\times \C^{N'})\not \subset \Sigma_1'$ for every integer $r$. Without loss of generality, we may assume that this subsequence is the whole sequence $(f^\ell)_\ell$. Furthermore, in what follows, when complexifying germs at $0\in \R^N$ of real-analytic mappings, we will keep the same notation for the obtained germs at $0\in \C^N$ of holomorphic maps.  For every integer $\ell$, we thus have 
\begin{equation}\label{e:dgrst}
 \rho_1(Z,\zeta,f^\ell(Z),\bar f^\ell (\zeta))\not \equiv 0,\ (Z,\zeta)\in (\6M,0).
\end{equation}
Since for every integer $\ell$, the mapping $\Lambda^k=T^\ell(u)$, $\xi'=g^\ell(x,u)$, $\omega'=h^\ell(x,u)$ satisfies the system of complex-valued real-algebraic equations given by \eqref{e:mainsystem}, we have, in view of Lemma~\ref{l:systT} that for every $\nu=0,\ldots,p_{\rho_1}-1$ that
\begin{equation}\label{e:invent}
 \6Q^{\rho_1}_\nu(Z,\zeta,T^\ell(u),\overline{T^\ell}(u),\xi',\varpi')\big|_{\6M \times \C^{2(N'-e)}}= 0,\ {\rm near}\ \,  0,
\end{equation}
where $\6Q^{\rho_1}$ is the polynomial associated to $\rho_1$ and given by Lemma~\ref{l:miami} and written as in \eqref{e:neweq}. Note also that we also have for $j=1,\ldots,N'-e$ and for $Z\in \C^N$ sufficiently close to the origin
$$\6P_j(Z,T^\ell(u),\overline{T^\ell}(u),g^\ell(Z),h^\ell_j(Z))=0.$$ 
Observe that the sequence of germs of holomorphic mappings $(f^\ell(Z))_\ell$ and $(T^\ell(u))_\ell$ converge in the Krull topology as $\ell \to \infty$ to the holomorphic maps $f(Z)=(g(Z),h(Z))$ and $(j^kf)(0,u)$ respectively. This implies in particular that, in view of \eqref{e:assure}, that for $\ell$ large enough, we have
$$ \6P_{j,m_j}\left(Z,T^\ell(u),\overline{T^\ell}(u),g^\ell(Z)\right)\not \equiv 0,\ j=0,\ldots,N'-e.$$
 Hence by Lemma~\ref{l:miami}, for $\ell$ large enough and for all $(Z,\zeta)\in \C^{2N}$ sufficiently to the origin, we have
\begin{equation}\label{e:calcul}
\6Q^{\rho_1} (Z,\zeta,T^\ell(u),\overline{T^\ell}(u),g^\ell(Z),\overline{g^\ell}(\zeta),\rho_1 (Z,\zeta,f^\ell(Z),\overline{f^\ell}(\zeta)))=0.
\end{equation}
Then \eqref{e:calcul} and \eqref{e:invent} imply that, for $(Z,\zeta)\in \6M$ near $0$, we have
\begin{equation}\label{e:arrival}
 \sum_{\nu=p_{\rho_1}}^{\delta}\6Q^{\rho_1}_\nu(Z,\zeta,T^\ell(u),\overline{T^\ell}(u),g^\ell(Z),\overline{g^\ell}(\zeta))\left(\rho_1(Z,\zeta,f^\ell (Z),\overline{f^\ell}(\zeta)) \right)^\nu= 0.
\end{equation}
Using \eqref{e:dgrst}, we may therefore write for $(Z,\zeta)\in \6M$ near $0$,
\begin{multline}\label{e:cnrst}
 \6Q^{\rho_1}_{p_{\rho_1}}(Z,\zeta,T^\ell(u),\overline{T^\ell}(u),g^\ell(Z),\overline{g^\ell}(\zeta))=\\ \sum_{\nu=1+p_{\rho_1}}^{\delta}\6Q^{\rho_1}_\nu(Z,\zeta,T^\ell(u),\overline{T^\ell}(u),g^\ell(Z),\overline{g^\ell}(\zeta))\left(\rho_1(Z,\zeta,f^\ell (Z),\overline{f^\ell}(\zeta)) \right)^\nu.
\end{multline}
Since ${\rm Graph}\, f \cap (M\times \C^{N'})\subset \Sigma'$ $($as germs at $(0,f(0)))$, we have $\rho_1(Z,\zeta,f(Z),\bar{f}(\zeta))\equiv 0$ for  $(Z,\zeta)\in \6M$ near $0$. Hence the restriction to $\6M$ of the right hand side of \eqref{e:cnrst} converges in the Krull topology, as $\ell \to \infty$, to $0$, whereas the left hand side of \eqref{e:cnrst} converges to $$\6Q^{\rho_1}_{p_{\rho_1}}(Z,\zeta,(j^kf)(0,u),(j^k \bar f)(0,u),g(Z),\bar{g}(\zeta))\big |_{\6M}.$$
This implies that  for $(Z,\zeta)\in \6M$ sufficiently close to $0$, we have
\begin{equation}\label{e:maroc}
\6Q^{\rho_1}_{p_{\rho_1}}(Z,\zeta,(j^kf)(0,u),(j^k \bar f)(0,u),g(Z),\bar{g}(\zeta))= 0.
\end{equation}
By the definition of $p_{\rho_1}$ (see \eqref{e:ictp}), we have  $$\6Q^{\rho_1}_{p_{\rho_1}}(Z,\zeta,(j^kf)(0,u),(j^k\bar f)(0,u),\xi',\varpi')\big|_{\6M \times \C^{2e}} \not \equiv 0,\quad {\rm near}\ 0.$$
Proposition~\ref{p:buffet} implies that the components of $g$ are algebraically dependent over the field $\4S^f$. This contradicts the definition of the mapping $g$. The proof of Proposition~\ref{p:summarize} is complete.
\end{proof}

\section{Proof of Theorem~\ref{t:nageneralbis}}\label{s:final}

Consider first the case $N=1$. Then $M$ is a real-algebraic curve in the complex plane. Fix a point $p\in M$. Since we can algebraically flatten $M$ near $p$, we may assume without loss of generality that $M$ is a piece of the real line and that $p=0$. Then the conclusion of the theorem follows easily from Theorem~\ref{t:artin} with $\4K=\R$.

Assume in the remainder of this section that $N\geq 2$. We first prove Theorem~\ref{t:nageneralbis} in the case where $M$ is generic in $\C^N$. Pick a point $p\in M$ and choose normal coordinates $Z=(z,w,u)\in \C^n \times \C^{d-c}\times \C^c$ vanishing at $p$ as in Lemma~\ref{l:coordconstantorbit}. As in Section~\ref{s:exam}, let $e$ be the transcendence degree of the field extension $\4S^f\hookrightarrow \4S^f(f)$ where $\4S^f$ is the quotient field of the ring $\4B^f$ as given in Definition~\ref{d:bf}. We use in what follows the notation introduced in the previous section. Let $S'$ be the real-algebraic subset of $\C^N \times \C^{N'}$ containing the germ at $(0,f(0))$ of ${\rm Graph}\, f \cap (M\times \C^{N'})$.

{\sc First case.} If $e=N'$, we claim that $M\times \C^{N'}\subset S'$ (as germs at $(0,f(0))$). Indeed, choose a finite collection of real-valued polynomials $r_1,\ldots,r_m$ over $\C^N \times \C^{N'}$ such that $S'=\{r_1=0,\ldots,r_m=0\}$. Suppose that there exists $j\in \{1,\ldots,m\}$ such that $r_j(Z,\bar Z,Z',\bar Z')\big|_{M\times \C^{N'}}\not \equiv 0$ near the origin. Since $r_j(Z,\zeta,f(Z),f(\zeta))\equiv 0$ for $(Z,\zeta)\in (\6M,0)$, Proposition~\ref{p:buffet} implies that the components of the mapping $f$ are algebraically dependent over $\4S^f$, which contradicts the fact that $e=N'$.  Hence for every integer $j\in \{1,\ldots,m\}$, $r_j(Z,\bar Z,Z',\bar Z') \equiv 0$ for $(Z,Z')\in (M\times \C^{N'},0)$, which proves the claim.
Therefore, we see that, if for every integer $\ell$, $f^\ell$ denotes the Taylor polynomial of the mapping $f$ up to order $\ell$, then the sequence of polynomial mappings $(f^\ell)_\ell$ satisfies the required conclusion.

{\sc Second case.} Assume that $e<N'$. In view of Proposition~\ref{p:summarize}, we see that it is enough to find for every integer $\ell$, germs at $0\in \R^N$ of complex-valued real-analytic mappings $\Lambda^k=T^\ell(u)$, $\xi'=g^\ell(x,u)$, $\omega'=h^\ell(x,u)$  that solve the system \eqref{e:mainsystem} with the following additional two properties: the sequences $(T^\ell(u))_\ell$, $(g^\ell(x,u))_\ell$, $(h^\ell(x,u))_\ell$ converge as $\ell \to \infty$ in the Krull topology to $(j^kf)(0,u)$, $g(x,u)$ and $h(x,u)$ respectively and the real and imaginary part of every component of these mappings must belong to the ring $\6A_{\R}\{x,u\}$. But the existence of these three sequences of germs of real-analytic mappings with the desired properties follows by a direct application of Theorem~\ref{t:popescu} (with $\4K=\R$) to the system of real-valued real-algebraic equations associated to the system \eqref{e:mainsystem}. The proof of Theorem~\ref{t:nageneralbis} is therefore complete in case $M$ is generic in $\C^N$.

To conclude the proof of Theorem~\ref{t:nageneralbis}, we will show as in \cite{MMZ3} that the case where $M$ is not generic follows from the generic case treated above. Indeed, if $M$ is not generic and $p\in M$, then the germ $(M,p)$ is equivalent (through a local complex-algebraic biholomorphism) to a germ of real-algebraic submanifold of the form $(M_1\times \{0\},0)\subset \C_t^{N-r}\times \C_s^r$ where $r\in \{1,\ldots,N-1\}$, $M_1$ is a real-algebraic generic submanifold in $\C^{N-r}$ of constant orbit dimension near $0$(see e.g.\ \cite{BERbook}). We may therefore assume that $(M,p)=(M_1\times \{0\},0)$ and apply the generic case treated above to the holomorphic map $(\C^{N-r},0)\ni t\mapsto f(t,0)\in \C^{N'}$. Hence, for every integer $\ell$, there exists a (germ of a) complex-algebraic mapping $\psi^\ell \colon (\C^{N-r},0)\to \C^{N'}$ such that ${\rm Graph}\, \psi^\ell\cap (M_1\times \C^{N'})\subset S'$ and that agrees with the power series mapping $f(t,0)$ up to order $\ell$ at $0$. Denoting, for every integer $\ell$, by $\varphi^\ell$ the Taylor polynomial of order $\ell$ of the mapping $(t,s)\mapsto f(t,s)-f(t,0)$ and setting $f^\ell:=\psi^\ell+\varphi^\ell$, we get the desired result.


\begin{thebibliography}{MMZ02a}


\bibitem[A69]{A69}
   M.~Artin: Algebraic approximation of structures over complete local rings. {\em Inst. Hautes \'Etudes Sci. Publ. Math.} {\bf 36} (1969), 23--58.



\bibitem[BER96]{BER96}
   M.S.~Baouendi; P.~Ebenfelt; L.P.~Rothschild: Algebraicity of holomorphic mappings between real algebraic sets in $\C^n$. {\em Acta Math.} {\bf 177} (1996), 225--273.




\bibitem[BER99]{BERbook}
   M.S.~Baouendi; P.~Ebenfelt; L.P.~Rothschild: 
  {\em Real Submanifolds in Complex Space and
  Their Mappings}. Princeton Math. Series {\bf 47}, 
  Princeton Univ. Press, 1999.





\bibitem[BMR02]{BMR}
   M.S.~Baouendi; N. Mir; L.P.~Rothschild: Reflection ideals and mappings between generic submanifolds in complex space. {\em J. Geom. Anal.} {\bf 12} (2002), 543--580.



\bibitem[BRZ01a]{BRZ1}
   M.S.~Baouendi; L.P.~Rothschild; D. Zaitsev: Equivalences of real submanifolds in complex space. {\em J. Differential Geom.} {\bf 59} (2001), 301--351.

\bibitem[BRZ01b]{BRZ2}
   M.S.~Baouendi; L.P.~Rothschild; D. Zaitsev: Points in general position in real-analytic submanifolds in $\C^N$ and applications.  {\em Complex analysis and geometry (Columbus, OH, 1999)},  1--20, Ohio State Univ. Math. Res. Inst. Publ., 9, de Gruyter, Berlin, 2001.


\bibitem[BCH08]{BCH08}
   S.~Berhanu; P.~Cordaro; J.~Hounie: {\em An introduction to involutive structures}. New Mathematical Monographs, {\bf 6}. Cambridge University Press, Cambridge, 2008.




\bibitem[CMS99]{CMS99}
   B.~Coupet; F.~Meylan; A.~Sukhov: Holomorphic maps of algebraic CR manifolds. {\em Int. Math. Res. Not.}  {\bf 1}, (1999), 1--29.

\bibitem[H94]{H94} X. Huang : On the mapping problem for algebraic real hypersurfaces in the complex spaces of different dimensions. {\em Ann. Inst. Fourier (Grenoble)}  {\bf 44}(2), (1994), 433--463.


\bibitem[HJ98]{HJ98} X.~Huang; S.~Ji : Global holomorphic extension of a local map and a Riemann mapping theorem for algebraic domains. {\em Math. Res. Lett.}  {\bf 5}, (1998), 247--260.


\bibitem[LM09]{LM09}
   B. Lamel; N. Mir: Holomorphic versus algebraic equivalence for deformations of real-algebraic CR manifolds,  {\em Comm. Anal. Geom.}, (31 pages, to appear).

\bibitem[MMZ03a]{MMZ3}
   F. Meylan; N. Mir; D. Zaitsev: Approximation and convergence of formal CR-mappings. {\em Int. Math. Res. Not.} {\bf 4} (2003), 211--242.


\bibitem[MMZ03b]{MMZsurvey}
   F. Meylan; N. Mir; D. Zaitsev: On some rigidity properties of mappings between CR-submanifolds in complex space.  {\em Journ\'ees ``\'Equations aux D\'eriv\'ees Partielles''},  Exp. No. XII, 20 pp., {\em Univ. Nantes, Nantes}, 2003. .



\bibitem[M98]{M98}
   N. Mir: Germs of holomorphic mappings between real-algebraic hypersurfaces, {\em Ann. Inst. Fourier (Grenoble)} {\bf 48}(4) (1998), 1025--1043.



\bibitem[P86]{P}
   D.~Popescu: General N\'eron desingularization and approximation, {\em Nagoya Math. J.} {\bf 104} (1986), 85--115.

\bibitem[W77]{W77} S.M. Webster: On the mapping problem for algebraic real hypersurfaces. {\em Invent. Math.} {\bf 43} (1977), 53-68.

\bibitem[Z99]{Z99} D.~Zaitsev:  Algebraicity of local holomorphisms between real-algebraic
submanifolds of complex spaces. {\em Acta Math.} {\bf 183} (1999),
273--305.



\bibitem[ZS58]{ZS58}
   O.~Zariski; P.~Samuel: 
  {\em Commutative Algebra, Volume I}. The University Series in Higher Mathematics, Springer-Verlag, New York, 1958.



\end{thebibliography}
\end{document}